\documentclass[aos,MSNbibl,seceqn,nameyear,dvips]{arximspdf}
\usepackage{mathbh}
\usepackage{graphicx}

\doi{10.1214/14-AOS1237} %kopijuoti is PTS
\volume{42}
\issue{4}
\pubyear{2014}
\firstpage{1598}
\lastpage{1634}
\docsubty{FLA}

\makeatletter
\newcommand{\rrvert}{\vert}
\newcommand{\llvert}{\vert}
\newtheorem{theorem}{Theorem}[section]
\newtheorem{lemma}[theorem]{Lemma}
\newproclaim{defin}[theorem]{Definition}
\newtheorem{corollary}[theorem]{Corollary}
\newtheorem{proposition}[theorem]{Proposition}
\newproclaim{example}[theorem]{Example}
\newtheorem{addendum}[theorem]{Addendum}
\newproclaim{condition}[theorem]{Condition}
\makeatother

\begin{document}
\begin{frontmatter}

\title{When uniform weak convergence fails: Empirical processes for
dependence functions and residuals~via epi- and hypographs\thanksref{T1}}
\runtitle{Weak convergence via epi- and hypographs}

\begin{aug}
\author[a]{\fnms{Axel}~\snm{B\"ucher}\corref{}\ead[label=e1]{axel.buecher@rub.de}},
\author[b]{\fnms{Johan}~\snm{Segers}\ead[label=e2]{johan.segers@uclouvain.be}}
\and
\author[a]{\fnms{Stanislav}~\snm{Volgushev}\ead[label=e3]{stanislav.volgushev@rub.de}}
\runauthor{A. B\"ucher, J. Segers and S. Volgushev}
\affiliation{Ruhr-Universit\"at Bochum, Universit\'e catholique de
Louvain and Ruhr-Universit\"at Bochum}
\address[a]{A. B\"ucher\\
S. Volgushev\\
Fakult\"at f\"ur Mathematik\\
Ruhr-Universit\"at Bochum\\
Universit\"atsstra{\ss}e 150\\
44780 Bochum\\
Germany\\
\printead{e1}\\
\phantom{E-mail: }\printead*{e3}}%adresu isvedimo komanda gale!
\address[b]{J. Segers\\
Institut de statistique\\
Universit\'e catholique de Louvain\\
Voie du Roman Pays, 20\\
B-1348 Louvain-la-Neuve\\
Belgium\\
\printead{e2}}
\end{aug}
\thankstext{T1}{Supported by the collaborative research center
``Statistical modeling of nonlinear dynamic processes'' (SFB 823)
of the German Research Foundation (DFG), by IAP research network Grant
No.~P7/06 of the Belgian government
(Belgian Science Policy) and by contract ``Projet d'Actions de
Recherche Concert\'ees'' No.~12/17-045
of the ``Communaut\'e fran\c{c}aise de Belgique.''}

% HISTORY:
\received{\smonth{5} \syear{2013}}
\revised{\smonth{2} \syear{2014}}

% ABSTRACT
%
\begin{abstract}
In the past decades, weak convergence theory for stochastic processes
has become a standard
tool for analyzing the asymptotic properties of various statistics.
Routinely, weak convergence is considered in the space of bounded
functions equipped with the supremum metric. However, there are cases
when weak convergence in those spaces fails to hold.
Examples include empirical copula and tail dependence processes and residual
empirical processes in linear regression models in case the underlying
distributions lack a certain degree of smoothness. To resolve the
issue, a new metric for locally bounded functions is introduced and the
corresponding weak convergence theory is developed. Convergence with
respect to the new metric is related to epi- and hypo-convergence and is
weaker than uniform convergence. Still, for continuous limits, it is
equivalent to locally uniform convergence, whereas under mild side
conditions, it implies $L^p$ convergence. For the examples mentioned
above, weak convergence with respect to the new metric is established
in situations where it does not occur with respect to the supremum distance.
The results are applied to obtain asymptotic properties of resampling
procedures and goodness-of-fit tests.\looseness=-1
\end{abstract}

% KEYWORDS
% Pirmas kwd is didziosios raides
%
\begin{keyword}[class=AMS]
\kwd[Primary ]{60F05}
\kwd{62G30}
\kwd[; secondary ]{62G32}
\kwd{62M09}
\end{keyword}
\begin{keyword}
\kwd{Bootstrap}
\kwd{copula}
\kwd{epigraph}
\kwd{hypograph}
\kwd{linear regression}
\kwd{local alternative}
\kwd{power curve}
\kwd{residual empirical process}
\kwd{stable tail dependence function}
\kwd{weak convergence}
\end{keyword}
\end{frontmatter}

%s1 #&#
\section{Introduction}

The Hoffman--J{\o}rgensen weak convergence theory in the space of
bounded functions is a great success story in mathematical
statistics\vadjust{\goodbreak}
[\citet{vandwell1996,kosorok2008}]. Measurability assumptions are
reduced to a minimum, no smoothness assumptions on the trajectories are
needed, it applies in a vast variety of circumstances, and the topology
of uniform convergence is fine enough so that, through the continuous
mapping theorem and functional delta method, it implies weak
convergence of a countless list of interesting statistical functionals.

But precisely because of the strength of uniform convergence, there are
circumstances where it does not hold. Weak convergence can fail when
the (pointwise) candidate limit process has discontinuous trajectories.
%%\sv{not precise, but perhaps the referee will understand this better}
%These arise when the candidate limit process (in the sense of
%convergence of finite-dimensional distributions) has discontinuous
%trajectories.
%For uniform convergence to take place, the locations of the
%discontinuities must be matched exactly. For the empirical
%distribution function based on a random sample from a discrete
%distribution, this holds. However, in other cases, it does not.
Think of empirical distributions based on residuals of some sort, that
is, observations that are themselves approximations of some latent
random variables. Because of the measurement error in the ordinates,
jump locations fail to be located exactly, and uniform convergence fails.
The examples which interest us in this paper concern empirical copula
processes, the empirical process based on residuals in a linear
regression setting and empirical tail dependence function processes.

A radical solution to the lack-of-convergence issue is to seek for
another metric on an appropriate function space. The metric should be
weak enough so that convergence does take place, but still strong
enough to enable statistical applications. Ideally, when the limit
process has continuous trajectories, it should be equivalent to uniform
convergence, so that in standard situations, nothing is lost. This is a
difficult task, and it turns out that the various extensions of
Skorohod's metrics [\citet{skorohod1956}] to functions of
several variables [\citet
{neuhaus1971,bickelwichura1971,straf1972,basspyke1985}] are not
suitable for the examples that we consider.

In the present paper, we construct such a metric by building on ideas
that originate in variational analysis and optimization theory.
% There, two modes of convergence, tailored towards the analysis of
%maximization and minimization problems, have been introduced.
In the context of minimization problems, one identifies a real function
$f$ on a suitable metric space $\mathbb{T}$ with its \emph{epigraph},
which is the set of all points $(x, y)$ in $\mathbb{T}\times\mathbb
{R}$ such that $f(x) \le y$.
% The epigraph is closed if and only if the function is lower
%semicontinuous.
\emph{Epi-convergence} of functions is then defined as Painlev\'e--Kuratowski convergence of their epigraphs [\citet
{beer1993,rockwets1998,molchanov2005}]. For maximization problems,
\emph{hypographs} and \emph{hypo-convergence} are defined in the same
way, the inequality sign pointing in the other direction.

Combining these modes of convergence, we will say that $f_n$ \emph
{hypi-converges} to $f$ if the epigraphs of $f_n$ converge to the
closure of the epigraph of $f$ and the hypographs of $f_n$ converge to
the closure of the hypograph of $f$. This mode of convergence is to be
distinguished from epi/hypo-convergence, a concept arising in
connection with saddle points [\citet{attouchwets1983}].

Broadly speaking, hypi-convergence is intermediate between uniform
convergence and $L^p$ convergence. Hypi-convergence implies uniform
convergence on compact subsets of the domain that are contained in the
set of continuity points of the limit function.
%of the domain at every point of which the limit function is continuous.
Hence, for continuous limits, we are back to uniform convergence. But
even without continuity, hypi-convergence implies convergence of global
extrema. Moreover, for limit functions which are continuous almost
everywhere, hypi-convergence implies $L^p$ convergence on compact sets.

In a similar way as one can consider weak epi-convergence of random
lower semicontinuous functions [\citet{geyer1994,molchanov2005}], we develop Hoffman--J{\o}rgensen weak convergence
theory with respect to the hypi-(semi)metric.
Thanks to an extension of the continuous mapping theorem for \mbox{semimetric}
spaces, we are able to leverage the above properties of
hypi-convergen\-ce to yield weak convergence of finite-dimensional
distributions, Kolmogorov--Smirnov type statistics, and procedures
related to $L^p$ spaces, notably Cram\'er--von Mises statistics. An
extension of the functional delta method is also discussed.

We investigate weak convergence with respect to the hypi-semimetric of
empirical copula processes, empirical tail dependence functions, and
empirical processes based on regression residuals. Weak
hypi-convergence of the empirical copula process is established for
copulas whose partial derivatives exist and are continuous everywhere
except for an arbitrary Lebesgue null set of the unit cube, which
extends results from \citet{rueschendorf1976} and others (see
Section~\ref{seccop} for more references on the empirical copula process).
% We apply weak hypi-convergence in the context of empirical copula
%processes and empirical tail dependence functions. For copulas whose
%partial derivatives exist and are continuous everywhere except for an
%arbitrary Lebesgue null set of the unit cube, weak hypi-convergence of
%the empirical copula process is established.
%This contains the theory of weak convergence in the supremum distance
%under stronger smoothness conditions as a special case.
From there, we show validity of the bootstrap [see \citet
{ferradweg2004}] and extend results on power curves for goodness-of-fit
tests under local alternatives [\citet{genquerem2007}]. Similar
results are shown for tail dependence functions, extending \citet
{buecdett2011} and \citet{einkraseg2012}. Classical results on
the empirical distribution function of regression residuals
[\citet{koul1969,loynes1980}] are extended to the case where the
true distribution has a discontinuous density.

The structure of the paper is as follows. The hypi-topology is
introduced in Section~\ref{secbase}. Weak convergence in hypi-space
is the topic of Section~\ref{secweak}. We provide tools for
checking weak hypi-convergence
%, including a variant of Slutsky's lemma,
and for exploiting it in a statistical context. The new framework is
applied for empirical copula processes in Section~\ref{seccop}, for
empirical tail dependence function processes in Section~\ref
{sectail}, and for the empirical process of regression residuals in
Section~\ref{secreg}. These three sections can be read independently
of one another.
%Section~\ref{secconcl} concludes.
In order to preserve the flow of the text, a number of auxiliary
results and all proofs are deferred to a sequence of \hyperref[APPVERIF]{Appendices} and an
online supplement [\citet{BucSegVol14supp}]. The weak convergence
theory for semimetric spaces in Appendix~\ref{APPWEAKSEMI}, including
a version of the extended continuous mapping theorem and the functional
delta method, is perhaps of independent interest.\looseness=-1

%%%%%%%%%%%%%%%%%%%%%%%%%%%%%%%%%%%%%%%%
%%% HYPI CONVERGENCE
%%%%%%%%%%%%%%%%%%%%%%%%%%%%%%%%%%%%%%%%
%s2 #&#
\section{Hypi-convergence of locally bounded functions}\label{secbase}

We introduce a mode of convergence for real-valued, locally bounded
functions on a locally compact, separable metric space (Section~\ref
{subsechypi}). For continuous limits, the metric is equivalent to
locally uniform convergence, but for discontinuous limits, it is
strictly weaker, while still implying $L^p$ convergence (Section~\ref
{subsecconseq}). The proofs for the results in this section are given
in Appendix~F.1.

%s2.1 #&#
\subsection{The hypi-semimetric}\label{subsechypi}

Let $(\mathbb{T}, d)$ be a locally compact, separable metric space.
% In our applications, $\TT$ will be a convex subset of Euclidean space.
The space $\mathbb{T}\times\mathbb{R}$ is a locally compact,
separable metric space, too, when equipped, for instance, with the metric
$d_{\mathbb{T}\times\mathbb{R}}((x_1, y_1), (x_2, y_2)) = \max\{
d(x_1, x_2), |y_1 - y_2| \}$.

Let $\ell_{\mathrm{loc}}^\infty(\mathbb{T})$ denote the space of
locally bounded functions $f\dvtx  \mathbb{T}\to\mathbb{R}$, that is,
functions that are uniformly bounded on compacta. If $\mathbb{T}$ is
itself compact, we will simply write $\ell^\infty(\mathbb{T})$.
Functions $f \in\ell_{\mathrm{loc}}^\infty(\mathbb{T})$ will be
identified with subsets of $\mathbb{T}\times\mathbb{R}$ by
considering their epigraphs and hypographs:
\begin{eqnarray*}
\operatorname{epi}f &=& \bigl\{ (x, y) \in\mathbb{T}\times\mathbb{R}\dvtx
f(x) \le y
\bigr\},
\\
\operatorname{hypo}f &=& \bigl\{ (x, y) \in\mathbb{T}\times\mathbb{R}\dvtx  y
\le f(x)
\bigr\}.
\end{eqnarray*}
Except for being locally bounded, functions $f$ in $\ell_{\mathrm{loc}}^\infty(\mathbb{T})$ can be arbitrarily rough. A minimal amount
of regularity will come from the lower and upper semicontinuous hulls
${f_\wedge} \le f \le{f_\vee} $:
%e2.1 #&#
%e2.2 #&#
%
\begin{eqnarray}
\label{eqfinf} {f_\wedge} (x) &=& \sup_{\varepsilon> 0} \inf\bigl
\{ f\bigl(x'\bigr)\dvtx  d\bigl(x', x\bigr) < \varepsilon
\bigr\},
\\
\label{eqfsup} {f_\vee} (x) &=& \inf_{\varepsilon> 0} \sup\bigl
\{ f\bigl(x'\bigr)\dvtx  d\bigl(x', x\bigr) < \varepsilon
\bigr\},
\end{eqnarray}
functions which are elements of $\ell_{\mathrm{loc}}^\infty(\mathbb
{T})$, too. Note that $(-f)_\wedge= -f_\vee$. A~convenient link
between epi- and hypographs on the one hand and lower and upper
semicontinuous hulls on the other hand is that
\begin{eqnarray*}
\operatorname{cl}(\operatorname{epi}f) &=& \operatorname{epi}
{f_\wedge},\qquad \operatorname{cl}(\operatorname{hypo}f) =
\operatorname{hypo} {f_\vee},
\end{eqnarray*}
where ``$\operatorname{cl}$'' denotes topological closure, in this
case, in the space $\mathbb{T}\times\mathbb{R}$. In particular, a
function $f$ is lower (upper) semicontinuous if and only if its
epigraph (hypograph) is closed.

Functions being identified with sets, notions of set convergence can be
applied to define convergence of functions. We rely on classical theory
exposed in, among others, \citet{matheron1975}, \citet
{beer1993}, \citet{rockwets1998} and \citet{molchanov2005}.
A standard topology on the space of closed subsets of a topological
space is the Fell hit-and-miss topology. If the underlying space is
locally compact and separable, as in our case, then the Fell topology
is metrizable. Moreover, in that case, convergence of a sequence of
closed sets in the Fell topology is equivalent to their Painlev\'
e--Kuratowski convergence. Recall that (not necessarily closed) sets
$A_n$ of a topological space converge to a set $A$ in the Painlev\'
e--Kuratowski sense if and only if (i) for every $x \in A$ there exists
a sequence $x_n$ with $x_n \in A_n$ such that $x_n \to x$ and
(ii)\vadjust{\goodbreak}
whenever $x_{n_k} \in A_{n_k}$ for some subsequence $n_k$ converges to
a limit $x$, we must have $x \in A$. The limit set $A$ is necessarily
closed, and Painlev\'e--Kuratowski convergence of $A_n$ to $A$ is
equivalent to Painlev\'e--Kuratowski convergence of $\operatorname
{cl}(A_n)$ to $A$.
% If the underlying space is compact, then Painlev\'e--Kuratowski
%convergence of closed sets is equivalent to their convergence in the
%Hausdorff metric.

Let $\mathcal{F}(\mathbb{T}\times\mathbb{R})$ be the space of
closed subsets of $\mathbb{T}\times\mathbb{R}$ and let ${d_{\mathcal
{F}}} $ be a metric inducing the Fell topology, or equivalently,
Painlev\'e--Kuratowski convergence. Examples of metrics ${d_{\mathcal
{F}}} $ for the Fell topology are to be found in \citet
{rockwets1998}, \citet{molchanov2005} and \citet
{ogura2007}. A~versatile notion of convergence of functions in
optimization theory is epi-convergence: a sequence of functions $f_n\dvtx
\mathbb{T}\to\mathbb{R}$ is said to epi-converge to a function $f$
if and only if the Painlev\'e--Kuratowski limit of $\operatorname
{epi}f_n$ (or equivalently, its closure) is equal to $\operatorname
{epi}f$, that is, if $d_{\mathcal{F}}( \operatorname
{cl}(\operatorname{epi}f_n), \operatorname{cl}( \operatorname{epi}f)
) \to0$ as $n \to\infty$. Necessarily, the limit set $\operatorname
{epi}f$ is closed and, therefore, $f$ must be lower semicontinuous.
Similarly, hypo-convergence of functions is defined as Painlev\'
e--Kuratowski convergence of their hypographs (or their closures), and
the limit function is necessarily upper semicontinuous.

If $f_n$ both epi-converges to ${f_\wedge} $ and hypo-converges to
${f_\vee} $, then we say that $f_n$ \emph{hypi-converges} to $f$.
This mode of convergence is the one that we propose in this paper.
According to the following result, hypi-convergence is metrizable and
can be checked conveniently by pointwise criteria.
% A formal definition is as follows.

%pr2.1 #&#
%
\begin{proposition}[(Hypi-convergence)]\label{prophypiconv}
Let $f_n, f \in\ell_{\mathrm{loc}}^\infty(\mathbb{T})$. The
following statements are equivalent:
\begin{longlist}[(iii)]
\item[(i)] $f_n$ epi-converges to ${f_\wedge} $ and hypo-converges to~${f_\vee} $.

\item[(ii)] The following pointwise criteria hold:
%e2.3 #&#
%
\begin{equation}
\label{eqhypiconvepi} \cases{ \forall x \in\mathbb{T}\dvtx  \forall x_n
\to x\dvtx
\displaystyle{f_\wedge} (x) \le\liminf_{n \to\infty}
f_n(x_n),
\vspace*{5pt}\cr
\forall x \in\mathbb{T}\dvtx  \exists
x_n \to x\dvtx  \displaystyle\limsup_{n \to\infty}
f_n(x_n) \le{f_\wedge} (x)}
\end{equation}
and
%e2.4 #&#
%
\begin{equation}
\label{eqhypiconvhypo} \cases{\forall x \in\mathbb{T}\dvtx  \forall x_n
\to x\dvtx
\displaystyle\limsup_{n \to\infty} f_n(x_n)
\le{f_\vee} (x),
\vspace*{5pt}\cr
\forall x \in\mathbb{T}\dvtx  \exists x_n
\to x\dvtx  \displaystyle{f_\vee} (x) \le\liminf_{n \to\infty}
f_n(x_n).}
\end{equation}
%
% \[
% \forall x \in\TT\dvtx  \forall x_n \to x\dvtx
% \finf(x) \le
% \liminf_{n \to\infty} f_n(x_n) \le
% \limsup_{n \to\infty} f_n(x_n) \le
% \fsup(x),
% \]
% and
% \[
% \begin{array}{l}
% \forall x \in\TT\dvtx  \exists x_n \to x\dvtx
% \displaystyle\limsup_{n \to\infty} f_n(x_n) \le\finf(x), \\[1em]
% \forall x \in\TT\dvtx  \exists x_n \to x\dvtx
% \displaystyle\fsup(x) \le\liminf_{n \to\infty} f_n(x_n).
% \end{array}
% \]

\item[(iii)]
The distance ${d_{\mathrm{hypi}}} (f_n, f)$ converges to $0$, where
${d_{\mathrm{hypi}}} $ denotes the \emph{hypi-semimetric} defined as
\[
{d_{\mathrm{hypi}}} (f, g) = \max\bigl\{ {d_{\mathcal{F}}} (
\operatorname{epi}
{f_\wedge}, \operatorname{epi} {g_\wedge} ),
{d_{\mathcal{F}}} ( \operatorname{hypo} {f_\vee},
\operatorname{hypo} {g_\vee} ) \bigr\},
\]
and ${d_{\mathcal{F}}} $ is a metric on $\mathcal{F}(\mathbb
{T}\times\mathbb{R})$ inducing the Fell topology.

\item[(iv)]
$f_n$ converges to $f$ in the \emph{hypi-topology}, which is defined
as the coarsest topology on $\ell_{\mathrm{loc}}^\infty(\mathbb
{T})$ for which the map
\[
\label{eqembedding} \ell_{\mathrm{loc}}^\infty(\mathbb{T}) \to\mathcal{F}(
\mathbb{T}\times\mathbb{R}) \times\mathcal{F}(\mathbb{T}\times\mathbb
{R})\dvtx  f
\mapsto\bigl( \operatorname{cl}(\operatorname{epi}f), \operatorname{cl}(
\operatorname{hypo}f) \bigr)
\]
is continuous,
% where $\mathcal{F}(\TT\times\reals)$ is equipped with the Fell
%topology and with the product topology on the Cartesian product,
that is, the hypi-open sets in $\ell_{\mathrm{loc}}^\infty(\mathbb
{T})$ are the inverse images of open sets in $\mathcal{F}(\mathbb
{T}\times\mathbb{R}) \times\mathcal{F}(\mathbb{T}\times\mathbb{R})$.
\end{longlist}
%
% In this case, we say that $f_n$ \emph{hypi-converges} to $f$.
\end{proposition}

Note that in (\ref{eqhypiconvepi}) and (\ref{eqhypiconvhypo}), we
can replace $f_n$ by $f_{n,\wedge}$ and $f_{n,\vee}$, respectively
(Lemma~\ref{lemconthull}). The equivalence of (i) and (ii) follows
from well-known pointwise criteria for epi- and hypo-convergence
[\citet{molchanov2005}, Chapter~5, Proposition~3.2(ii)].
Statements (iii) and (iv) are just reformulations of what it means to
have both epi- and hypo-convergence in (i).

%f1 #&#
%
\begin{figure}%[ht!]

\includegraphics{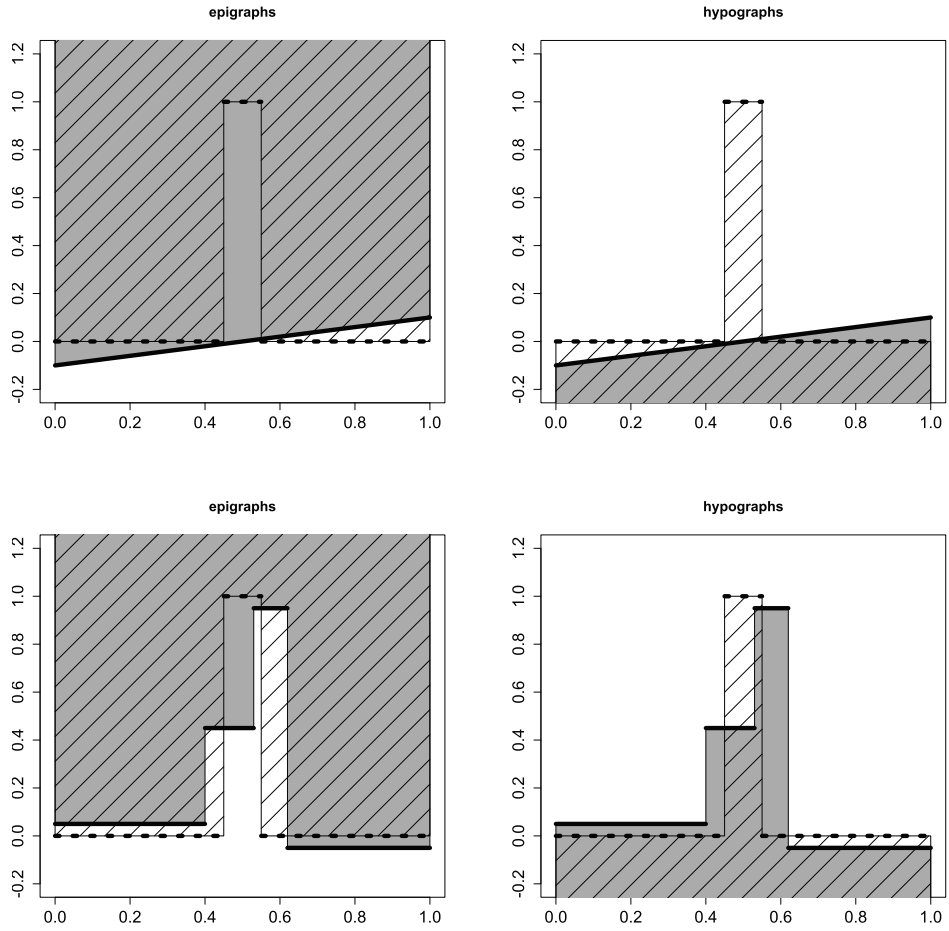}

\caption{\textup{Top}: two functions whose epigraphs are close but whose
hypographs are far away.
\mbox{\textup{Bottom}}: two functions of which both the epi-
and hypographs are close.
The light gray areas represent epigraphs (left column) and hypographs
(right column) of the functions depicted by solid lines, whereas shaded
areas represent epigraphs (left column) and hypographs (right column)
of the functions depicted by dotted lines.}\label{fig2}
\end{figure}

Intuitively, two functions are close in the hypi-semimetric if both
their epigraphs and their hypographs are close.
% Note that the key, and the novelty, of our approach is to
%simultaneously consider closeness of both epi- \textit{and}
%hypo-graphs.
Two functions on $\mathbb{T}= [0,1]$ whose epigraphs are close but
whose hypographs are far away are depicted in the upper part of
Figure~\ref{fig2}: for instance, the point $(0.5, 1)$ belongs to the
hypograph of the dotted-line function but is far away from any point in
the hypograph of the solid-line function. As a consequence, these two
functions are not close in the hypi-semimetric. For comparison, two
functions that are close in the hypi-semimetric are depicted in the
lower part of Figure~\ref{fig2}.

By Proposition~\ref{prophypiconv}(ii), if $f_n$ hypi-converges to $f$
and if $f$ is continuous at $x$, then $f_n(x_n) \to f(x)$ whenever $x_n
\to x$.
%The pointwise criteria can be rephrased as follows\dvtx  for every $x \in
%hypo-converges at $x$ to $\fsup(x)$, in the sense of equations~
Moreover, it follows that locally uniform convergence of locally
bounded functions implies their hypi-convergence. The converse is not
true if the hypi-limit is not continuous.

Hypi-convergence of sequences $f_n$ and $g_n$ to $f$ and $g$,
respectively, does in general not imply hypi-convergence of the
sequence of sums, $f_n + g_n$, to the sum of the limits, $f + g$. For
instance, let $x_n$ converge to $x$ in $\mathbb{T}$ with $x_n \ne x$
and set $f_n = \mathbh{1}_{\{ x_n \}}$ and $g_n = - \mathbh{1}_{\{x\}
}$. Still, a sufficient condition is that at least one of the limit
functions is continuous; see Lemma~\ref{lemsum} for an even more
general result.

% For c\`adl\`ag functions on the interval $[0, 1]$, it can be seen
%that hypi-convergence is equivalent to Skorohod's $M_2$ convergence [
%completed graphs $\Gamma_f = \{(x, y)\dvtx  \finf(x) \le y \le\fsup(x)
%of $[0, 1] \times\RR$. Another way of viewing the latter convergence
%is as graphical convergence of set-valued mappings $x \mapsto[
%general domains, this characterization fails.

An alternative view on the hypi-topology can be gained by identifying
$f\in\ell_{\mathrm{loc}}^\infty(\mathbb{T})$ with its completed
graph $\Gamma(f) = \operatorname{epi}({f_\wedge} ) \cap
\operatorname{hypo}({f_\vee} )$ [\citet{vervaat1981}]. We
suspect that for certain domains, hypi-convergence is equivalent to set
convergence of completed graphs. For c\`adl\`ag functions on $\mathbb
{T}= [0,1]$, the latter convergence can be seen to be equivalent to
Skorohod $M_{2}$-convergence [\citet{molchanov2005}, page~377],
whence hypi-convergence can be regarded as a coordinate-free extension
of Skorohod $M_{2}$-convergence to nonsmooth functions on rather
general domains.

%s2.2 #&#
\subsection{Leveraging hypi-convergence}
\label{subsecconseq}

As mentioned already, uniform convergence implies hypi-convergence but
not conversely. Nevertheless, at subsets of the domain where the limit
function is continuous, the converse does hold. In this sense, working
in hypi-space does not necessarily yield weaker results than in the
uniform topology. All proofs for this section are given in
Appendix~F.1 in the supplement [\citet{BucSegVol14supp}].

%pr2.2 #&#
%
\begin{proposition}
\label{propK}
Let $K \subset\mathbb{T}$ be compact and let $f \in\ell_{\mathrm{loc}}^\infty(\mathbb{T})$ be continuous at every $x \in K$. If $f_n$
hypi-converges to $f$ in $\ell_{\mathrm{loc}}^\infty(\mathbb{T})$,
then $\sup_{x \in K} |f_n(x) - f(x)| \to0$ as $n \to\infty$.
\end{proposition}

Being a combination of epi- and hypo-convergence, hypi-convergence
preserves convergence of extrema. Later, we will make use of this
property when investigating Kolmogorov--Smirnov type test statistics
(Section~\ref{subsecpower}).

%pr2.3 #&#
%
\begin{proposition}
\label{propextrema}
Let $G\subset\mathbb{T}$ be an open subset with compact closure. If
$f_n$ hypi-converges to $f$ in $\ell_{\mathrm{loc}}^\infty(\mathbb
{T})$ and if $f$ is continuous on the boundary of $G$, then $\inf
f_n(G) \to\inf f(G)$ and $\sup f_n(G) \to\sup f(G)$
as $n \to\infty$. If $G=\mathbb{T}$ is compact, then the boundary of
$G$ in $\mathbb{T}$ is empty, and hence the conclusions hold true
without imposing any conditions on $f$.
%Let $\TT$ be compact. If $f_n$ hypi-converges to $f$ in $\ell^\infty(
\end{proposition}

Hypi-convergence implies $L^p$-convergence for finite $p$, provided
that the limit function is not too rough. This is useful, for instance,
when studying Cram\'er--von Mises statistics (Section~\ref
{subsecpower}) and other statistical procedures based on the
$L^2$-distance, such as minimum distance estimators. Recall that upper
and lower semicontinuous functions are necessarily Borel measurable.

%pr2.4 #&#
%
\begin{proposition}
\label{propLp}
Let $\mu$ be a finite Borel measure supported on a compact subset of
$\mathbb{T}$. If $f_n$ hypi-converges to $f$ in $\ell^\infty(\mathbb
{T})$ and if $f$ is continuous $\mu$-almost everywhere, then, for
every $p \in[1, \infty)$, we have $\int|f_{n,\vee} - f_{n,\wedge
}|^p \,{d}\mu\to0$ and $\int|f_{n}^*-f^*|^p \,{d}\mu\to0$ as $n \to\infty$, where $f_{n}^*$ and $f^*$ represent
arbitrary Borel measurable functions on $\mathbb{T}$ such that
$f_{n,\wedge} \le f_n^* \le f_{n,\vee}$ and ${f_\wedge} \le f^* \le
{f_\vee} $.
% In particular, if both $f$ and each $f_n$ are Borel measurable, we
%have $\int|f_{n}-f|^p \diff\mu\to0$.
\end{proposition}

%%%%%%%%%%%%%%%%%%%%%%%%%%%%%%%%%%%%%%%%
%%% WEAK HYPI CONVERGENCE
%%%%%%%%%%%%%%%%%%%%%%%%%%%%%%%%%%%%%%%%
%s3 #&#
\section{Weak hypi-convergence of stochastic processes}
\label{secweak}

% \subsection{Weak convergence in the hypi-topology}

When applying\break  Hoffman--J{\o}rgensen weak convergence theory, it is
customary to work in a metric space. However, ${d_{\mathrm{hypi}}} $
is a semimetric and not a metric: if functions $f, g \in\ell_{\mathrm{loc}}^\infty(\mathbb{T})$ share the same lower and upper
semicontinuous hulls, then ${d_{\mathrm{hypi}}} (f, g) = 0$ even if
$f$ and $g$ are different functions.

To obtain a metric space, we consider equivalences classes of functions
at hypi-distance zero. For $f \in\ell_{\mathrm{loc}}^\infty(\mathbb
{T})$, let $[f]$ be the\vspace*{1pt} set of all $g \in\ell_{\mathrm{loc}}^\infty
(\mathbb{T})$ such that ${d_{\mathrm{hypi}}} (f, g) = 0$. Let
$L_{\mathrm{loc}}^\infty(\mathbb{T})$ be the\vspace*{1pt} space of all such
equivalence classes. Then $L_{\mathrm{loc}}^\infty(\mathbb{T})$
becomes a metric space when equipped with the hypi-metric (abusing
notation) ${d_{\mathrm{hypi}}} ([f], [g]):= {d_{\mathrm{hypi}}} (f,
g)$. The map $[ \cdot]$ from $\ell_{\mathrm{loc}}^\infty(\mathbb
{T})$ into $L_{\mathrm{loc}}^\infty(\mathbb{T})$ sending $f$ to
$[f]$ is continuous and it sends open sets to open sets and closed sets
to closed sets.

Let $X_n$ and $X$ be maps from probability spaces $\Omega_n$ and
$\Omega$, respectively, into $\ell_{\mathrm{loc}}^\infty(\mathbb
{T})$. Assume that $X$ is hypi Borel measurable, that is, measurable
with respect to the $\sigma$-field generated\vspace*{1pt} by the hypi-open sets of
$\ell_{\mathrm{loc}}^\infty(\mathbb{T})$. Then the map $[X] = [
\cdot](X)$ into $L_{\mathrm{loc}}^\infty(\mathbb{T})$ is Borel
measurable, too. Since $L_{\mathrm{loc}}^\infty(\mathbb{T})$ is a
metric space, weak convergence theory as in \citet{vandwell1996}
applies: we say that $X_n$ weakly hypi-converges to $X$ in $\ell
_{\mathrm{loc}}^\infty(\mathbb{T})$ if and\vspace*{1pt} only if $[X_n]
\rightsquigarrow[X]$ in $L_{\mathrm{loc}}^\infty(\mathbb{T})$.
Simplifying notation, we sometimes omit brackets and write $X_n
\rightsquigarrow X$ in $L_{\mathrm{loc}}^\infty(\mathbb{T})$.

In order to prove weak hypi-convergence, we will usually combine an
initial result on weak convergence of some stochastic process, usually
some empirical process and with respect to the supremum distance, with
the (extended) continuous mapping theorem [\citet{vandwell1996},
Theorems~1.3.6 and~1.11.1]. The task then consists of proving
hypi-cont\-inuity of the relevant mappings into $\ell_{\mathrm{
loc}}^\infty(\mathbb{T})$ on sufficiently large subsets of their
domains. Two situations of particular importance are the following:
\begin{itemize}
\item
convergence of sums, inducing in particular a variant of Slutsky's
lemma (Lemma~\ref{lemSlutsky} and Appendix~\ref{appsum});
\item
convergence to a function in $\ell_{\mathrm{loc}}^\infty(\mathbb
{T})$ that is defined as the upper or lower semicontinuous hull of some
other function that is originally defined on a dense subset of $\mathbb
{T}$ only (Appendix~\ref{appext}).
\end{itemize}

% The continuous mapping theorem for semimetric spaces (Theorem~
%weak hypi-convergence in $\Loc(\TT)$.

%le3.1 #&#
%
\begin{lemma}[(Slutsky)]
\label{lemSlutsky}
%Let $(\TT,\dT)$ be a locally compact, separable metric space.
Let $X_n, Y_n\dvtx  \Omega_n \to\ell_{\mathrm{loc}}^\infty(\mathbb
{T})$ be arbitrary maps and let $X\dvtx  \Omega\to\ell_{\mathrm{loc}}^\infty(\mathbb{T})$ be Borel measurable with respect to the
hypi-semimetric. If $[X_n] \rightsquigarrow[X]$ and $[Y_n]
\rightsquigarrow[0]$ in $L_{\mathrm{loc}}^\infty(\mathbb{T})$, then
$[X_n + Y_n] \rightsquigarrow[X]$ in $L_{\mathrm{loc}}^\infty
(\mathbb{T})$.
\end{lemma}

The proof of this and all other results from this section are given in
Appendix~F.2 in the supplement [\citet{BucSegVol14supp}].

By Proposition~\ref{propK}, the map from $(\ell_{\mathrm{loc}}^\infty(\mathbb{T}), {d_{\mathrm{hypi}}} )$ into $(\ell
^\infty(K), \| \cdot\|_\infty)$ sending a function $f$ to its
restriction $f|_K$ on a compact subset $K$ of $\mathbb{T}$ is
continuous at every function $f$ which is continuous in every point $x$
of $K$; here $\| \cdot\|_\infty$ denotes the supremum norm. By a
generalization of the continuous mapping theorem to semimetric spaces
(Theorem~\ref{thmcmt}), weak hypi-convergence implies weak
convergence with respect to the supremum distance insofar the limit
process has continuous trajectories. More precisely, we have the
following result.

%co3.2 #&#
%
\begin{corollary}
\label{corK}
Let $X_n$ and $X$ be maps from probability spaces $\Omega_n$ and
$\Omega$, respectively, into $\ell_{\mathrm{loc}}^\infty(\mathbb
{T})$ such that $X$ is hypi Borel measurable. If $[X_n]
\rightsquigarrow[X]$ in $L_{\mathrm{loc}}^\infty(\mathbb{T})$ and
if $K \subset\mathbb{T}$ is a nonempty, compact set such that, with
probability one, $X$ is continuous in every $x \in K$, then $X_n|_K
\rightsquigarrow X|_K$ in $(\ell^\infty(K), \| \cdot\|_\infty)$.
\end{corollary}

Taking $K$ to be finite, we find that weak hypi-convergence implies
weak convergence of finite-dimensional distributions at points where
the limit process is continuous almost surely.

For finite Borel measures $\mu$ on $\mathbb{T}$ with compact support,
Proposition~\ref{propLp} states the continuity of the embedding from
the set of Borel measurable functions of $\ell_{\mathrm{loc}}^\infty
(\mathbb{T})$ equipped with the hypi-topology into $L^p(\mu)$, for
every $1 \le p < \infty$. Again by the continuous mapping theorem
(Theorem~\ref{thmcmt}), weak hypi-convergence then implies weak
$L^p$-convergence. A technical nuisance is that in order to view
$L^p(\mu)$ as a metric space, we have to consider equivalence classes
of functions that are equal $\mu$-almost everywhere; notation $[
\cdot]_\mu$.

%co3.3 #&#
%
\begin{corollary}
\label{corLp}
Let $X_n$ and $X$ be maps from probability spaces $\Omega_n$ and
$\Omega$, respectively, into $\ell_{\mathrm{loc}}^\infty(\mathbb
{T})$ such that $X$ is hypi Borel measurable. Let~$\mu$ be a finite
Borel measure on $\mathbb{T}$ with compact support. If $[X_n]
\rightsquigarrow[X]$ in $L_{\mathrm{loc}}^\infty(\mathbb{T})$ and
if $X$ is $\mu$-almost everywhere continuous with probability one,
then $\int|X_{n,\vee} - X_{n,\wedge}|^p \,{d}\mu$
converges to $0$ in outer probability and both $[X_{n,\vee}]_\mu$ and
$[X_{n,\wedge}]_\mu$ converge weakly in $L^p(\mu)$ to $[X_\vee]_\mu
= [X_\wedge]_\mu$, for every $p \in[1, \infty)$.
\end{corollary}

% \begin{remark}
% Addition not being continuous on $\loc(\TT)$, the latter space is not
%a topological vector space. This prohibits a direct application of the
%functional delta method to weak hypi-convergence. In practice, this
%difficulty can be circumvented by first mapping $\loc(\TT)$ or subsets
%thereof into a better-behaving space by means of the previous
%corollaries and then apply the functional delta method on the image
%space. An alternative approach is to work out an abstract version of
%the functional delta method to vector spaces equipped with a
%(semi)metric for which addition is not continuous everywhere.
% \end{remark}

Addition not being continuous on $\ell_{\mathrm{loc}}^\infty(\mathbb
{T})$, the latter space is not a topological vector space.
%Even worse, addition of equivalence classes in $\Loc(\TT)$ is not even
%well-defined, which
This prohibits a direct application of the functional\vadjust{\goodbreak} delta method
[\citet{vandwell1996}, Theorem~3.9.4] to weak hypi-convergence.
However, in Appendix~\ref{APPWEAKSEMI}, we provide a variant of the
functional delta method (Theorem~\ref{thmdelta}) that is sufficiently
flexible to deal with maps defined on semimetric spaces endowed with an
addition operator that is not necessarily continuous.

%certain empirical processes that do not necessarily converge with
%respect to the supremum distance.}

%%%%%%%%%%%%%%%%%%%%%%%%%%%%%%%%%%%%%%%%
%%% COPULAS
%%%%%%%%%%%%%%%%%%%%%%%%%%%%%%%%%%%%%%%%
%s4 #&#
\section{Empirical copula processes}\label{seccop}

Usually, empirical copula processes are studied in the space of bounded
functions on $[0, 1]^d$ equipped with the supremum distance. Weak
convergence then requires existence and continuity of the first-order
partial derivatives of the copula on the interior and some subsets of
the boundary of $[0, 1]^d$. In this section, we show what can be done
in case the latter smoothness condition is not satisfied. Existence and
continuity of the partial derivatives almost everywhere is still enough
to ensure weak hypi-convergence of the empirical copula process
(Section~\ref{subseccop}). The result is strong enough to validate
the bootstrap (Section~\ref{subsecbootcop}) and to analyze
Kolmogorov--Smirnov and Cram\'er--von Mises test statistics, even under
local alternatives (Section~\ref{subsecpower}). The proofs of the
results in this section are given in Appendix~\ref{appcop} and,
partially, in Appendix~F.3 in the supplement
[\citet{BucSegVol14supp}].

%s4.1 #&#
\subsection{Weak hypi-convergence}
\label{subseccop}

Let $\mathbf{X}_i=(X_{i1},\ldots, X_{id})$, with $i \in\mathbb{N}$, be a
strictly stationary sequence of $d$-variate random vectors. (No
confusion should arise from the use of the symbol ``$d$'' for both the
metric on $\mathbb{T}$ above and the dimension of the random vectors
here.) Throughout this section, the joint distribution function $F$ of
$\mathbf{X}_i$ is assumed to have continuous marginal distributions $F_1,\ldots, F_d$ and its copula is denoted by $C$.
Further, for $j=1,\ldots, d$, let $U_{ij}=F_j(X_{ij})$
% be the probability integral transform of the components of $
and set $\mathbf{U}_i=(U_{i1},\ldots, U_{id})$. Note that $\mathbf{U}_i$ is
distributed according to~$C$. Consider the empirical distribution functions
\[
F_n(\mathbf{x})= \frac{1}{n}\sum_{i=1}^n
\mathbh{1}\{ \mathbf{X}_i\leq\mathbf{x}\}, \qquad G_n(\mathbf{u})=
\frac{1}{n}\sum_{i=1}^n \mathbh{1}\{
\mathbf{U}_i \leq\mathbf{u} \}
\]
for $\mathbf{x}\in\mathbb{R}^d$ and $\mathbf{u}\in[0,1]^d$. For a
distribution function $H$ on the reals, let
\begin{eqnarray*}
H^{-}(u):=\cases{ \displaystyle\inf\bigl\{x\in\mathbb{R}\dvtx  H(x)\geq u
\bigr\}, &\quad$0<u\leq1$,
\vspace*{5pt}\cr
\displaystyle\sup\bigl\{x\in\mathbb{R}\dvtx  H(x)=0\bigr
\}, &\quad$u=0$,}
\end{eqnarray*}
denote the (left-continuous) generalized inverse function of $H$.

The object of interest is the empirical copula, defined by
\[
C_n(\mathbf{u}) = F_n\bigl(F_{n1}^-(u_1),\ldots, F_{nd}^-(u_d)\bigr), \qquad\mathbf{u} \in[0,
1]^d,
\]
where $F_{nj}$ denotes the $j$th marginal empirical distribution
function. For convenience, we will abbreviate the notation for the
empirical copula by $C_n(\mathbf{u})= F_n(\mathbf F_{n}^-(\mathbf{u}) )$, with
$\mathbf F_n^-(\mathbf{u})=(F_{n1}^-(u_1),\ldots, F_{nd}^-(u_d))$.

Often, the empirical copula is defined as the distribution function of
the vector of rescaled ranks, and/or it is turned into a genuine copula
via linear interpolation. Since these variants often differ from the
empirical copula by at most a term of order $o_p(n^{-1/2})$, uniformly
over $[0, 1]^d$, they do not affect the asymptotic distribution of the
empirical copula process, defined by
%e4.1 #&#
%
\begin{equation}
\label{eqccn} \mathbb{C}_n = \sqrt n ( C_n - C).
\end{equation}
The\vspace*{1pt} asymptotic behavior of $\mathbb{C}_n$, especially its weak
convergence in the space $\ell^\infty([0, 1]^d)$ equipped with the
supremum norm $\| \cdot\|_\infty$, has been investigated by several
authors under various conditions [\citet{rueschendorf1976,ferradweg2004,ghoudiremillard2004,tsukahara2005,vandervaartwellner2007,deheuvels2009,segers2012,buecvolg2013}].

The main arguments to derive the limit of $\mathbb{C}_n$ are as
follows. For the sake of a clear explanation, let us assume for the
moment that the random vectors $(\mathbf{X}_i)_{i\in\mathbb{N}}$ form an
i.i.d. sequence, even though the same arguments work for many time
series models with short-range dependence. Observing that $C_n=F_n(\mathbf
F_n^-) = G_n (\mathbf G_n^-)$, we can decompose $\mathbb{C}_n$ into two terms:
%e4.2 #&#
%
\begin{equation}
\label{eqmd} \mathbb{C}_n = \sqrt n\bigl\{ G_n\bigl(
\mathbf G_n^-\bigr)- C \bigr\} = \alpha_n\bigl(\mathbf
G_n^-\bigr) + \sqrt n \bigl\{ C\bigl(\mathbf G_n^-\bigr) - C
\bigr\},
\end{equation}
where $\alpha_n=\sqrt n ( G_n - C)$ denotes the usual empirical
process associated to the sequence $(\mathbf{U}_i)_{i\in\mathbb{N}}$.

Deriving the limit of the first term in (\ref{eqmd}) is standard:
since $\alpha_n \rightsquigarrow\alpha$ in $\ell^\infty([0,1]^d)$
with respect to the supremum norm, for a $C$-Brownian bridge $\alpha$,
and since $\sup_{0 \le u_j \le1} | G_{nj}^-(u_j) - u_j| = o_p(1)$, we
obtain that $\alpha_n(\mathbf G_n^-) \rightsquigarrow\alpha$ in $(\ell
^\infty([0,1]^d), \| \cdot\|_\infty)$, too.

Regarding the second term in (\ref{eqmd}), the argumentation is
harder. Set $\beta_n=(\beta_{n1},\ldots, \beta_{nd})$, where $\beta
_{nj} = \sqrt n ( G_{nj}^- - \mathrm{id}_{[0,1]})$ denotes the
quantile process of the $j$th coordinate and where $\mathrm{id}_A$ is
the identity map on a set $A$. It follows from the functional delta
method applied to the inverse mapping $H \mapsto H^-$ that $\| \beta
_{nj}+\alpha_{nj} \|_\infty=o_p(1)$, where $\alpha_{nj}(u_j) =
\alpha_n(1,\ldots, 1, u_j, 1,\ldots, 1)$, with $u_j\in[0,1]$ at the
$j$th position. Therefore, $\beta_{nj} \rightsquigarrow-\alpha_{j}$
in $(\ell^\infty([0,1]), \| \cdot\|_\infty)$, where, similarly,
$\alpha_j(u_j)$ is defined as $\alpha(1,\ldots, 1, u_j, 1, \ldots,
1)$. Now, %we can write
%e4.3 #&#
%
\begin{equation}
\label{eq2ndterm} \sqrt n \bigl\{ C\bigl(\mathbf G_n^-\bigr) - C\bigr\} =
\sqrt n \bigl\{ C(\mathrm{id}_{[0,1]^d} + \beta_n / \sqrt n) - C
\bigr\},
\end{equation}
which can be handled under suitable differentiability conditions on
$C$. To conclude upon weak convergence with respect to the supremum
distance, the weakest assumption so far has been stated in \citet
{segers2012}.

%co4.1 #&#
%
\begin{condition}
\label{condpd}
For $j=1,\ldots,d$ the partial derivatives $\dot C_j(\mathbf{u})$ exist
and are continuous on $\{\mathbf{u} \in[0,1]^d\dvtx  u_j \in(0,1)\}$.
\end{condition}

Under Condition~\ref{condpd},
%e4.4 #&#
%
\begin{equation}
\label{eqdiffterm} \sqrt n \bigl\{ C\bigl(\mathbf G_n^-\bigr) - C\bigr\}
(\mathbf
u) \rightsquigarrow- \sum_{j=1}^d \dot
C_j(\mathbf u) \alpha_j(u_j),
\end{equation}
where $\dot C_j(\mathbf u)$ can be defined, for instance, as $0$ if $u_j \in
\{0, 1\}$. Hence,
%e4.5 #&#
%
\begin{equation}
\label{eqCsup} \mathbb{C}_n(\mathbf u) \rightsquigarrow\mathbb{C}(\mathbf u)
= \alpha(\mathbf u) - \sum_{j=1}^d \dot
C_j(\mathbf u) \alpha_j(u_j)
\end{equation}
in $\ell^\infty([0, 1]^d)$ with respect to the supremum distance.

Condition~\ref{condpd} ensures that the limit process $\mathbb{C}$
in (\ref{eqCsup}) has continuous trajectories. Actually, if $\mathbb
{C}_n$ is to converge weakly with respect to the supremum distance,
then the weak limit must have continuous trajectories with probability
one. The reason is that the mapping
\[
\Delta\dvtx  \ell^\infty\bigl([0, 1]^d\bigr) \to[0, \infty)\dvtx  f
\mapsto\sup_{\mathbf u \in[0, 1]^d} \bigl| {f_\vee} (\mathbf u) -
{f_\wedge} (\mathbf u) \bigr|
\]
is continuous with respect to $\| \cdot\|_\infty$ and that $0 \le
\Delta\mathbb{C}_n \le d/\sqrt{n} \to0$ almost surely. The
expression for $\mathbb{C}$ in (\ref{eqCsup}) then suggests that
$\mathbb{C}_n$ does not converge weakly in $(\ell^\infty([0, 1]), \|
\cdot\|_\infty)$ if Condition~\ref{condpd} does not hold.

% to some procecontinuity of the trajectories of Cn converges to C
%supremum then limit is continuous
%
% If the condition is not satisfied, weak convergence with respect to
%the supremum distance can fail.

%ex4.2 #&#
%
\begin{example}[(Mixture model)]
\label{excop}
For $\lambda\in(0, 1)$, consider the bivariate copula given by
\[
C(u_1, u_2) = (1-\lambda) u_1
u_2 + \lambda\min(u_1, u_2).
\]
For $u_1 \ne u_2$, the partial derivatives are
\begin{eqnarray*}
\dot{C}_1(u_1, u_2) &=& (1-\lambda)
u_2 + \lambda\mathbh{1}(u_1 < u_2),
\\
\dot{C}_2(u_1, u_2) &=& (1-\lambda)
u_1 + \lambda\mathbh{1}(u_2 < u_1).
\end{eqnarray*}
On the diagonal $u_1 = u_2$, the partial derivatives do not exist.
Still, by the decomposition in (\ref{eqmd}), the finite-dimensional
distributions of $\mathbb{C}_n$ can be seen to converge to the ones of
the process $\tilde{\mathbb{C}}$ defined as
\[
\tilde{\mathbb{C}}(u_1, u_2) = \alpha(u_1,
u_2) - \dot{C}_1(u_1, u_2)
\alpha_1(u_1) - \dot{C}_2(u_1,
u_2) \alpha_2(u_2),
\]
if $u_1 \ne u_2$, whereas, on the diagonal $u_1 = u_2 = u$,
\[
\tilde{\mathbb{C}}(u, u) = \alpha(u, u) - (1-\lambda) u \bigl\{
\alpha_1(u) + \alpha_2(u) \bigr\} - \lambda\max\bigl(
\alpha_1(u), \alpha_2(u) \bigr),
\]
the distribution of which is non-Gaussian.

Now suppose that $\mathbb{C}_n \rightsquigarrow\mathbb{C}$ in $(\ell
^\infty([0, 1]^d, \| \cdot\|_\infty)$ for some $\mathbb{C}$. Then
the finite-dimensional distributions of $\mathbb{C}$ must be equal to
the ones of $\tilde{\mathbb{C}}$. Additionally, the trajectories of
$\mathbb{C}$ must be continuous almost surely, and thus the law of the
random variable $\mathbb{C}(u_1, u_2)$ must depend continuously on the
coordinates $(u_1, u_2)$. However, by the above expressions for $\tilde
{\mathbb{C}}$, continuity cannot hold at points on the diagonal. This
yields a contradiction and, therefore, $\mathbb{C}_n$ cannot converge
weakly with respect to the supremum distance.
\end{example}

By considering weak hypi-convergence, we can go far beyond
Condition~\ref{condpd}. Condition~\ref{condD} imposes the
regularity needed to deal with the left-hand side of (\ref
{eq2ndterm}) in the hypi-semimetric.

%co4.3 #&#
%
\begin{condition}
\label{condD}
The set $\mathcal S$ of points in $[0,1]^d$ where the partial
derivatives of the copula $C$ exist and are continuous has Lebesgue measure~$1$.
\end{condition}

Since a copula is monotone in each of its arguments, its partial
derivatives automatically exist almost everywhere. Condition~\ref
{condD} then only concerns continuity of these partial derivatives. In
practice, Condition~\ref{condD} poses no restriction at all.
Still, there do exist copulas that do not satisfy Condition~\ref
{condD}. It can be shown that a bivariate example is given by the
copula with Lebesgue density
\[
c(u,v) = \frac{\mathbh{1}_{A\times A}(u,v)}{\lambda_1(A)} + \frac
{\mathbh{1}_{B\times B}(u,v)}{\lambda_1(B)},
\]
where $\lambda_1$ denotes the one-dimensional Lebesgue measure, $A
\subset[0,1]$ is a closed set which is at the same time nowhere dense
and satisfies $\lambda_1(A) \in(0,1)$ and where $B = [0,1]\setminus A$.
%Still, there do exist copulas that do not satisfy Condition~
%conditional distribution of $U_2$ given $U_1 = u_1$ depend on whether
%$u_1$ belongs to a set $A \subset[0,1]$ or not, where $A$ is at the
%same time nowhere dense and has positive Lebesgue measure.

% This condition, which obviously implies Condition~\ref{condpd}, is
%already enough to prove weak hypi-convergence of the emprical copula
%process.
For broad applicability, we relax the assumption of serial independence
and replace it by the following condition, which holds for i.i.d.
sequences as well as for stationary sequences under various weak
dependence conditions [\citet{rio2000,doukhanfermanianlang2009,dehlingdurieu2011}].

%co4.4 #&#
%
\begin{condition} \label{condweak}
The empirical process $\alpha_n=\sqrt n (G_n - C)$ converges weakly in
$(\ell^\infty([0,1]^d), \|\cdot\|_\infty)$ to some limit process
$\alpha$ which has continuous sample paths, almost surely.
\end{condition}

Under Condition~\ref{condD}, the term on the right-hand side of (\ref
{eqdiffterm}) is defined only on $\mathcal{S}$. We extend it to the
whole of $[0, 1]^d$ by taking lower semicontinuous hulls as in
Appendix~\ref{appext}. Let $| \cdot|$ denote the Euclidean norm
and let $\mathcal{C}(A)$ be the set of continuous real-valued
functions on a domain $A$. Recall our convention of omitting the
brackets $[ \cdot ]$ when working in $L_{\mathrm{loc}}^\infty
(\mathbb{T})$.

%th4.5 #&#
%
\begin{theorem}
\label{theocopweak}
Suppose that Condition~\ref{condweak} holds and that $C$ satisfies
Condition~\ref{condD}. Then
%e4.6 #&#
%
\begin{equation}
\label{eqcc} \mathbb{C}_n \rightsquigarrow\mathbb{C}= \alpha+
{d}C_{(-\alpha_{1},\ldots,-\alpha_{d})} %\quad\mbox{ in }~ (L^
\end{equation}
in $(L^\infty([0,1]^d), {{d_{\mathrm{hypi}}} })$, where, for $a=(a_1,\ldots, a_d)\in\{ \mathcal{C}([0,1]) \}^d$,
\[
{d}C_a(\mathbf u ) =\sup_{\varepsilon>0} \inf\Biggl\{
\sum_{j=1}^d \dot{C}_j(\mathbf v)
a_j(v_j)\dvtx  \mathbf v \in\mathcal{S}, | \mathbf v - \mathbf u | <
\varepsilon\Biggr\}.
\]
\end{theorem}

By Section~\ref{secbase}, Theorem~\ref{theocopweak} has several
useful consequences.
\begin{itemize}
\item
First, it implies weak convergence with respect to the supremum
distance of the restriction of the empirical copula process to compact
subsets of the union of~$\mathcal{S}$ and the boundary of $[0, 1]^d$;
see Corollary~\ref{corK}. This is akin to the convergence results for
multilinear empirical copulas for count data in \citet
{genestneslehovaremillard2013}.
%The limit is the same as in \eqref{eqCsup}, and in particular we
%retrieve the usual weak convergence result under Condition~
Note that, in particular, we obtain the weak convergence result in
(\ref{eqCsup}) under the stronger Condition~\ref{condpd}.
% In particular, it implies the weak convergence result in
%Proposition~3.1 in \citet{segers2012}.
%
\item
Furthermore, we obtain weak convergence of the empirical copula process
in $(L^p([0,1]^d), \| \cdot\|_p)$ for any $1 \le p < \infty$. To the
best of our knowledge, this result is new and opens the door to
$L^p$-type inference procedures for a broad class of copulas.
\end{itemize}
Two possible applications are treated in the following subsections.

% The list of further possible applications of hypi-convergence of the
%empirical copula process is large, among them being for instance the
%derivation of the asymptotics of minimum distance estimators or
%goodness-of-fit tests. In the subsequent sections, we stick to two
%important illustrations in order to demonstrate the power of the new
%methodology.

%s4.2 #&#
\subsection{A bootstrap device}\label{subsecbootcop}

Assume that $\mathbf X_1,\ldots, \mathbf X_n$ are serially independent. We
show that the bootstrap based on resampling with replacement
[\citet{ferradweg2004}] and the bootstrap based on the multiplier
central limit theorem [\citet{buecherdette2010}] provide valid
approximations for $\mathbb{C}$ with respect to the hypi-semimetric.
Our multiplier bootstrap is different from the approach in \citet
{remillardscaillet2009}, which requires estimation of the first-order
partial derivatives of $C$.

Let $M\in\mathbb{N}$ be some large integer and, for each $m\in\{1,\ldots, M\}$, let
$\mathbf X_1^{[m]},\ldots,\break \mathbf X_n^{[m]}$ be drawn with
replacement from the sample. The \emph{resampling bootstrap empirical
copula process} is defined as
%e4.7 #&#
%
\begin{equation}
\label{eqccnm} \mathbb{C}_n^{[m]}= \sqrt{n}\bigl(C_n^{[m]}-
C_n\bigr),
\end{equation}
where $C_n^{[m]}$ denotes the empirical copula calculated from the
bootstrap sample $\mathbf X_1^{[m]},\ldots, \mathbf X_n^{[m]}$.
% Note that the distribution of $\CC_n^\tm$ is equal to the
%distribution of $\sqrt{n}(F_n^\tm(\vect F_n^{\tm-}) - C_n \}$,
Note that we can represent $C_n^{[m]}$ by $F_n^{[m]}(\mathbf F_n^{[m]-})$,
where
\[
F_n^{[m]}(\mathbf x) = \frac{1}{n} \sum
_{i=1}^n W_{ni}^{[m]}\mathbh{1}( \mathbf
X_i \le\mathbf x)
\]
and where $W_n^{[m]}=(W_{n1}^{[m]},\ldots, W_{nn}^{[m]})$ denotes a
multinomial random vector with $n$ trials, $n$ possible outcomes, and
success probabilities $(1/n,\ldots, 1/n)$, independent of the sample
and independent across $m \in\{1, \ldots, M\}$.

Regarding the multiplier bootstrap, let $\{\xi_{i}^{[m]}\dvtx  i \ge1, m
= 1, \ldots, M \}$ be i.i.d. random variables, independent of the
sample, with both mean and variance equal to one and such that $\int
_0^\infty\sqrt{ \mathbb{P}( \xi_i > x) } \,dx<\infty$. Let
\[
\tilde F_n^{[m]}(\mathbf x) = \frac{1}{n} \sum
_{i=1}^n \xi_{i}^{[m]}\mathbh{1}( \mathbf
X_i \le\mathbf x), \qquad\mathbf x \in\mathbb{R}^d
\]
and define
%e4.8 #&#
%
\begin{equation}
\label{eqccnmtilde} \tilde\mathbb{C}_n^{[m]} = \sqrt n \bigl\{ \tilde
F_n^{[m]}\bigl(\tilde{\mathbf{F}}_n^{[m]-}\bigr) -
C_n \bigr\}
\end{equation}
as the \emph{multiplier bootstrap empirical copula process}. The\vspace*{1pt}
following proposition shows that both $\mathbb
{C}_n^{[1]},\ldots, \mathbb
{C}_n^{[M]}$ and $\tilde\mathbb
{C}_n^{[1]},\ldots, \tilde\mathbb
{C}_n^{[M]}$ can be regarded as asymptotically
independent copies of $\mathbb{C}_n$.

%pr4.6 #&#
%
\begin{proposition}
\label{propbootweak}
Let $\mathbf{X}_i$, $i \in\mathbb{N}$, be i.i.d. $d$-variate random
vectors with common distribution function $F$ having\vspace*{1pt} continuous margins
and a copula $C$ satisfying Condition~\ref{condD}. Let $\mathbb
{C}_n$, $\mathbb{C}_n^{[m]}$ and $\tilde{\mathbb
{C}}_n^{[m]}$ be as in (\ref{eqccn}), (\ref
{eqccnm}) and (\ref{eqccnmtilde}), \mbox{respectively}. Then both
$(\mathbb{C}_n, \mathbb{C}_n^{[1]},\ldots,
\mathbb{C}_n^{[M]})$ and $(\mathbb{C}_n, \tilde
\mathbb{C}_n^{[1]},\ldots, \tilde\mathbb
{C}_n^{[M]})$\break  weakly\vspace*{1pt} converge to $(\mathbb{C},
\mathbb{C}^{[1]},\ldots, \mathbb
{C}^{[M]})$ in the space $(L^\infty([0,1]^d),
{{d_{\mathrm{hypi}}} })^{M+1}$, where $\mathbb
{C}^{[1]},\ldots, \mathbb
{C}^{[M]}$ denote independent copies of $\mathbb
{C}$ in (\ref{eqcc}).
% Under Condition~\ref{condD}, both
% $(\CC_n, \CC_n^{[1]},\ldots, \CC_n^{
%[M]})$ and $(\CC_n, \tilde\CC_n^{
%[1]},\ldots, \tilde\CC_n^{
%[M]})$ weakly converge to $(\CC, \CC^{[1]},\ldots,
\end{proposition}

By hypi-continuity of the supremum and infimum functionals (see
Proposition~\ref{propextrema}), the bootstrap approximation can, for
instance, be used to construct asymptotic uniform confidence bands for
the copula.

%s4.3 #&#
\subsection{Power curves of tests for independence}\label{subsecpower}

In the present section, we derive weak hypi-con\-vergence of the
empirical copula process for triangular arrays.
We apply it to the problem of comparing statistical tests for
independence by local power curves. This comparison has been carried
out by \citet{genquerem2007} under strong differentiability
assumptions on copula densities. By considering hypi-convergence, we
can extend their results to copulas that do not have a density with
respect to the Lebesgue measure.
% Analogous extensions are possible for comparing tests for different
%testing problems, see for instance \citet{bergques2009} for
%goodness-of-tests.

We consider a triangular array of random vectors $\mathbf X_1^{(n)},\ldots, \mathbf X_{n}^{(n)}$ which are row-wise i.i.d. with continuous
marginals and copula $C^{(n)}$. We suppose that there exists a copula
$C$ satisfying Condition~\ref{condD} such that
%e4.9 #&#
%
\begin{equation}
\label{eqloc} \Delta_n = \sqrt n \bigl\{ C^{(n)} - C \bigr
\} \to\Delta
\end{equation}
uniformly, for some continuous function $\Delta$ on $[0,1]^d$.
Let $C_n^{(n)}$ denote the empirical copula based on $\mathbf X_1^{(n)},\ldots, \mathbf X_{n}^{(n)}$. Let $\mathbf U_1^{(n)},\ldots, \mathbf U_{n}^{(n)
}$ denote the sample obtained by the marginal probability integral
transform and let $G_n^{(n)}$ and $\alpha_n^{(n)}$ denote its
empirical distribution function and empirical process,\vadjust{\goodbreak} respectively.
Similarly as before, we have the decomposition
\begin{eqnarray*}
\mathbb{C}_n^{(n)} &=& \sqrt n \bigl\{ C_n^{(n)}
- C^{(n)} \bigr\}
\\
&=&\sqrt n \bigl\{ G_n^{(n)}\bigl( \mathbf G_n^{(n)-}
\bigr) - C^{(n)}\bigl(\mathbf G_n^{(n)-}\bigr) \bigr\} +
\sqrt n \bigl\{ C^{(n)}\bigl( \mathbf G_n^{(n)-}\bigr) -
C^{(n)} \bigr\}
\\
&=& \alpha_n^{(n)} \bigl( \mathbf G_n^{(n)-}
\bigr) + \sqrt n \bigl\{ C\bigl( \mathbf G_n^{(n)-}\bigr) - C
\bigr\} + \bigl\{ \Delta_n\bigl( \mathbf G_n^{(n)-}
\bigr) - \Delta_n \bigr\}.
\end{eqnarray*}
We will\vspace*{1pt} show in Appendix~F.3 in the supplement [\citet
{BucSegVol14supp}] that $\alpha_n^{(n)} \rightsquigarrow\alpha$ in
$(\ell^\infty([0,1]^d), \|\cdot\|_\infty)$, where~$\alpha$ is a
$C$-Brownian bridge. Therefore,
%, and since $1 + 1 = 2$,
the first summand weakly converges to $\alpha$ with respect to the
supremum norm. The second summand weakly converges in the hypi-topology
to ${d}C_{(-\alpha_1,\ldots, -\alpha_d)}$, while the last
one converges to $\Delta-\Delta\equiv0$, uniformly. This motivates
the following result.

%pr4.7 #&#
%
\begin{proposition} \label{proploc}
Given the above set-up
%Under Condition~\ref{condweak}
and if (\ref{eqloc}) is met with $C$ satisfying Condition~\ref
{condD}, we have $ \mathbb{C}_n^{(n)} \rightsquigarrow\mathbb{C}$
in $(L^\infty([0,1]^d, {{d_{\mathrm{hypi}}} })$, where $\mathbb{C}$
is the same process as in Theorem~\ref{theocopweak}. Additionally, in
$(L^\infty([0,1]^d), {{d_{\mathrm{hypi}}} })$,
\[
\sqrt{n}\bigl( C_n^{(n)} - C\bigr) \rightsquigarrow
\mathbb{C}+\Delta.
\]
%
% in $(L^\infty([0,1]^d), {\dhypi})$.
\end{proposition}

To illustrate the latter result, we investigate the local efficiency of
tests for independence as considered in \citet{genquerem2007}.
Instead of imposing conditions (i) and (ii) on page~169 in their paper,
we only suppose that (\ref{eqloc}) holds with $C=\Pi$, the
independence copula, and $\Delta= \delta\Lambda$, where $\Lambda\in
\mathcal{C}([0,1]^d)$ and $\delta\ge0$.
For brevity, we only compare the test statistics
\[
T_n = n \int_{[0,1]^d} \bigl\{ C_n^{(n)}
- \Pi\bigr\}^2 \,d \Pi\quad\mbox{and}\quad S_n =
\sqrt n \bigl\| C_n^{(n)} - \Pi\bigr\|_\infty.
\]
From weak hypi-convergence of $\sqrt{n} (C_n^{(n)} - C)$ and
Propositions~\ref{propK} and \ref{propLp}, we obtain that
\[
T_n \rightsquigarrow\mathcal{T}_{\delta}= \int
_{[0,1]^d} ( \mathbb{C}+ \delta\Lambda) ^2 \,d
\Pi, \qquad S_n \rightsquigarrow\mathcal{S}_\delta= \|
\mathbb{C}+ \delta\Lambda\|_\infty.
\]
Hence, the local power curves of the tests to the level $\alpha\in
(0,1)$ in direction~$\Lambda$ are given by
\[
\delta\mapsto\mathbb{P}\bigl\{\mathcal{T}_\delta> q_{\mathcal
{T}_0}(1-
\alpha)\bigr\}, \qquad\delta\mapsto\mathbb{P}\bigl\{\mathcal{S}_\delta>
q_{\mathcal{S}_0}(1-\alpha)\bigr\},
\]
where $q_{\mathcal{T}_0}(1-\alpha)$ and $q_{\mathcal{S}_0}(1-\alpha
)$ denote the $(1-\alpha)$-quantiles of $\mathcal{T}_0$ and $\mathcal
{S}_0$, respectively. These curves can be compared by analytical
calculations as in \citet{genquerem2007} or by simulation.

\section{Stable tail dependence functions}\label{sectail}

Let $\mathbf{X}_1,\ldots, \mathbf{X}_n$, where $\mathbf{X}_i=\break  (X_{i1},\ldots,X_{id})$, be i.i.d.~$d$-variate random vectors with distribution
function $F$ and continuous marginal distribution functions $F_1,\ldots, F_d$. We assume that the following limit, called the stable tail
dependence function of $F$,
%e5.1 #&#
%
\begin{eqnarray}
\label{eqL}\quad  L(\mathbf{x}) &=& \lim_{t \downarrow0} t^{-1}
\mathbb{P} \bigl\{ 1- F_1(X_{11}) \le tx_1
\mbox{ or } \cdots\mbox{ or } 1-F_d(X_{1d}) \le
tx_d \bigr\},
\end{eqnarray}
exists as a function $L\dvtx [0,\infty)^d\rightarrow[0,\infty)$.

For $i\in\{1,\ldots, n\}$ and $j\in\{1,\ldots, d\}$, let $R_{i}^j$
denote the rank of $X_{ij}$ among $X_{1j},\ldots, X_{nj}$. Replacing
all distribution functions in (\ref{eqL}) by their empirical
counterparts and replacing $t$ by $k/n$ where $k=k_n$ is a positive
sequence such that $k_n\to\infty$ and $k_n=o(n)$, we obtain the
following nonparametric estimator for $L_n$, called the empirical
(stable) tail dependence function:
\[
\hat L_n (\mathbf{x}) = \frac{1}{k} \sum
_{i=1}^n \mathbh{1} \biggl\{ R_i^1
> n+\frac{1}{2}-kx_1 \mbox{ or } \cdots\mbox{ or }
R_i^d > n + \frac{1}{2} - kx_d/n
\biggr\}
\]
[\citet{huang1992,dreehuan1998}]. The inclusion of the term
$1/2$ inside the indicators serves to improve the finite sample
behavior of the estimator.

In \citet{einkraseg2012}, a functional central limit theorem for
$\sqrt{k} ( \hat{L}_n - L )$ is given in the topology of uniform
convergence on compact subsets of $[0, \infty)^d$. The result requires
$L$ to have continuous first-order partial derivatives on sufficiently
large subsets of $[0, \infty)^d$, similar to Condition~\ref{condpd}
for copulas. By switching to weak hypi-convergence, we are able to get
rid of smoothness conditions altogether.

Similarly as in Section~\ref{seccop}, let $\mathcal{S}$ denote the
set of all points $\mathbf x\in[0,\infty)^d$ where $L$ is differentiable.
\label{textS} The function $L$ being convex, Theorem~25.5 in
\citet{rockafellar1970} implies that the complement of $\mathcal
{S}$ is a Lebesgue null set and that the gradient $(\dot{L}_1,\ldots,
\dot L_d)$ of $L$ is continuous on $\mathcal{S}$. Proceeding as in
Appendix~\ref{appext}, we may define, for any $(a_1,\ldots, a_d) \in
\{C([0,\infty))\}^d$, a function on $[0, \infty)^d$ by
%e5.2 #&#
%
\begin{equation}
\label{eqdlw} {d}L_{(a_1,\ldots, a_d)}( \mathbf x ) = \sup
_{\varepsilon> 0}
\inf\Biggl\{ \sum_{j=1}^d
\dot{L}_j(\mathbf y) a_j(y_j)\dvtx  \mathbf y \in
\mathcal{S}, | \mathbf x - \mathbf y | < \varepsilon\Biggr\}.
\end{equation}

As in \citet{einkraseg2012}, let $\Lambda$ be the Borel measure
on $[0, \infty)^d$ such that $\Lambda( A(\mathbf x) ) = L(x)$ where
$A(\mathbf x) = \bigcup_{j=1}^d \{ y \in[0, \infty)^d\dvtx  y_j \le x_j \}$
for $\mathbf x \in[0, \infty)^d$. Let $\mathbb{W}$ be a mean-zero
Gaussian process on $[0, \infty)^d$ with continuous trajectories and
with covariance function $\mathbb{E}[ \mathbb{W}(\mathbf x) \mathbb
{W}(\mathbf y) ] = \Lambda( A(\mathbf x) \cap A( \mathbf y ) )$. Let $\Delta
_{d-1} = \{ \mathbf x \in[0, 1]^d\dvtx  x_1 + \cdots+ x_d = 1\}$ be the unit
simplex in $\mathbb{R}^d$. For $f\in\ell_{\mathrm{loc}}^\infty
([0,\infty)^d)$ and $j=1,\ldots,d$, define\vspace*{1pt} $f_j^{0}\in\ell_{\mathrm{loc}}^\infty([0,\infty))$ through $f_j^0(x_j)=f(0,\ldots,0,x_j,0,\ldots,0)$.
Recall our convention of omitting the brackets $[ \cdot ]$ when
working in $L_{\mathrm{loc}}^\infty(\mathbb{T})$.

%th5.1 #&#
%
\begin{theorem}
\label{theotailweak}
Let $\mathbf{X}_i$, $i \in\mathbb{N}$, be i.i.d. $d$-dimensional random
vectors with common distribution function $F$ with continuous margins
$F_1, \ldots, F_d$ and stable tail dependence function $L$. Suppose
that the following conditions hold:
\begin{longlist}[(C1)]
\item[(C1)] For some $\alpha>0$ we have, uniformly in $\mathbf{x} \in\Delta_{d-1}$,
\begin{eqnarray*}
&& t^{-1} \mathbb{P} \bigl\{ 1- F_1(X_{11}) \le
tx_1 \mbox{ or } \cdots\mbox{ or } 1-F_d(X_{1d})
\le tx_d \bigr\}
\\
&&\qquad = L(\mathbf x) + O\bigl(t^\alpha\bigr), \qquad t \downarrow0.
\end{eqnarray*}
\item[(C2)] We have $k=o(n^{2\alpha/(1+2\alpha)})$ and $k\rightarrow\infty
$ as $n\rightarrow\infty$.
\end{longlist}
Then, in $(L_{\mathrm{loc}}^\infty([0,\infty)^d), {d_{\mathrm{hypi}}} )$,
\[
\sqrt{k} (\hat L_n - L ) \rightsquigarrow\mathbb{W}+
{d}L_{(-\mathbb{W}_{1}^0,\ldots,-\mathbb{W}_{d}^0)}, \qquad n \to
\infty.
\]
\end{theorem}

The proof of Theorem~\ref{theotailweak} is similar to the one of
Theorem~\ref{theocopweak} and is deferred to the supplement
[\citet{BucSegVol14supp}].

Conditions~(C1) and (C2) also appear in Theorem~4.6 in \citet
{einkraseg2012} and are needed to ensure that the estimator is
asymptotically unbiased. The difference with their theorem is that we
do not need their condition~(C3) on the partial derivatives of $L$.
Therefore, Theorem~\ref{theotailweak} also covers piecewise linear
stable tail dependence functions arising from max-linear models
[\citet{wangstoev2011}].

Weak hypi-convergence of $\sqrt{k} ( \hat{L}_n - L )$ can be
exploited to validate statistical procedures for tail dependence
functions in the same way as was done with weak hypi-convergence of
empirical copula processes in Section~\ref{seccop}. In contrast to
copulas, no smoothness conditions on $L$ are needed at all.
Applications include the bootstrap [\citet{pengqi2008}] and
minimum $L^2$-distance estimation [\citet{buecdett2011}].
% In particular, Proposition~\ref{theminidicop} remains valid upon
%replacing $C$ by $L$ and $\sqrt{n}$ by $\sqrt{k}$.
Hypi-convergence implying $L^2$-convergence, Theorem~\ref
{theotailweak} also provides another way to prove the asymptotic
normality of the M-estimator in \citet{einkraseg2012}.

\section{Error distributions in regression models}\label{secreg}
% In this section, we
Consider a linear regression model for a sample $(\mathbf{X}_i, Y_i)$, $i
\in\{1, \ldots, n\}$, in $\mathbb{R}^p \times\mathbb{R}$, of the form
%e6.1 #&#
%
\begin{equation}
\label{eqregmod} Y_i = \mathbf{X}_i' \bolds{
\beta} + \varepsilon_i.
\end{equation}
Here, $(\mathbf{X}_i, \varepsilon_i)$, for $i \in\{1, \ldots, n\}$, are
i.i.d. random vectors in $\mathbb{R}^p \times\mathbb{R}$. It is
assumed that $\mathbf{X}_i$ and $\varepsilon_i$ are independent and that
the distribution of $\varepsilon_i$ is constrained in such a way that
the vector of regression coefficients $\bolds{\beta}$ is identifiable
(provided the support of $\mathbf{X}_i$ is sufficiently large). For
instance, the requirement $\mathbb{E}(\varepsilon_i) = 0$ yields a
mean regression model, whereas $\operatorname{median}(\varepsilon_i)
= 0$ yields a median regression model. For simplicity, we restrict
attention to serial independence and to a scalar dependent variable.

% The assumption that $(\vect{X}_i, \eps_i)_{i=1,...,n}$ is independent
%can be relaxed to incorporate time series models or general kinds of
%dependence. However, the basic phenomena which we describe in this
%section are going to remain the same, and we therefore state all
%results under this assumption. \sv{provide more details?} %Possible
%extensions will be indicated throughout the text. \sv{I would probably
%prefer to state generic assumptions}

The model is semiparametric with parametric component $\bolds{\beta} \in
\mathbb{R}^p$ and nonparametric components $P^X$ and $P^\varepsilon$,
the distributions of the explanatory variables $\mathbf{X}_i$ and the
errors $\varepsilon_i$. We are interested in the estimation of the
cumulative distribution function, $F$, of $\varepsilon_i$:
\[
F(z) = \mathbb{P}( \varepsilon_i \le z ), \qquad z \in\bar{
\mathbb{R}},
\]
where $\bar{\mathbb{R}} = [-\infty, \infty]$, a convenient
compactification of the real line.

Let $\hat{\bolds{\beta}}_n$ be a consistent estimator for $\bolds{\beta
}$. In Theorem~\ref{theoresweak} below, we will be more specific
about the asymptotic distribution of $\hat{\bolds{\beta}}_n$. We define
estimated residuals as
%e6.2 #&#
%
\begin{equation}
\label{eqhateps} \hat{\varepsilon}_{n,i} = Y_i -
\mathbf{X}_i' \hat{\bolds{\beta}}_n =
\varepsilon_i - \mathbf{X}_i' (\hat{\bolds{
\beta}}_n - \bolds{\beta})
\end{equation}
%
%Borrowing terminology from copula theory, we may coin these random
%variables `pseudo-observations' from $F$.
and obtain a simple estimator for $F$ by
%e6.3 #&#
%
\begin{equation}
\label{eqhatF} \hat{F}_n(z) = \frac{1}{n} \sum
_{i=1}^n \mathbh{1} (\hat{\varepsilon}_{n,i}
\le z), \qquad z \in\bar{\mathbb{R}}.
\end{equation}
The \emph{empirical residual process} corresponding to $\hat{F}_n$ is
%e6.4 #&#
%
\begin{equation}
\label{eqFn} \mathbb{F}_n(z) = \sqrt{n} \bigl\{
\hat{F}_n(z) - F(z) \bigr\}, \qquad z \in\bar{\mathbb{R}}.
\end{equation}
Weak convergence results for $\mathbb{F}_n$ play a central role in,
for example, testing the goodness-of-fit of error distributions or in
the derivation of the asymptotic behavior of more sophisticated
estimators for $F$; see \citet{koul2002} for an overview. First
results on the asymptotic behavior of $\mathbb{F}_n$ were derived in
\citet{koul1969} and \citet{loynes1980} (in generalized
regression models), and more recently those findings were extended in
various directions such as, for instance, time series analysis [see
\citet{koul2002,englernielsen2009}, and the references cited
therein] or coefficient vectors of growing dimension [see \citet
{chenlockhart2001}, for an overview].
%For the development of weak convergence results for $\FF_n$ in the
%space $\ell^\infty(\reals)$ \ab{$\loc(\reals)$} equipped with either
%the supremum distance or the hypi-semimetric, a crucial assumption
%will be that $F$ is absolutely continuous. Its probability density
%function, $f$, will show up in the expression of the limit process.
%In the literature hitherto, corresponding results were usually shown
%under the assumption that $F$ has a continuous probability density
%function $f$.
All of those extensions share the assumption that $F$ has a continuous
probability density function $f$. In that case, weak convergence takes
place with respect to the supremum distance
%process.
and the process admits an expansion of the form
%e6.5 #&#
%
\begin{equation}
\label{eqreprffn}\quad  \mathbb{F}_n(z) = \Biggl[ \frac{1}{\sqrt{n}}\sum
_{i=1}^n \bigl\{\mathbh{1}(
\varepsilon_i \le z) - F(z)\bigr\} \Biggr] + f(z) \mathbb{E}[
\mathbf{X}]' \sqrt{n} ( \hat{\bolds{\beta}}_n - \bolds{\beta} ) +
o_p(1)\hspace*{-20pt}
\end{equation}
uniformly in $z \in\bar{\mathbb{R}}$, where $\mathbf{X}$ denotes a
random vector with the same distribution as $\mathbf{X}_i$. In the present
section, we will drop the assumption of continuity of $f$ and consider
weak hypi-convergence of $\mathbb{F}_n$.

%The estimator $\hat{F}_n(z)$ may be improved by incorporating the
%identifiability constraint on the distribution of $\eps_i$. In a mean
%regression model, for instance, the knowledge that $\expec(\eps_i) =
%0$ may be exploited via empirical likelihood to construct a more
%accurate estimator. (If the regressors include an intercept term and
%if $\vect{\beta}$ is estimated via ordinary least squares, then the
%residuals already have empirial mean zero and such an improvement is
%not possible. How about other estimators?) However, the asymptotic
%theory for such more sophisticated estimators will usually be based on
%the one for the naive estimator $\hat{F}_n(z)$. \js{References: To be
%completed.}

%The estimator $\hat{F}_n(z)$ differs from the empirical distribution
%function, $F_n(z)$, based on the unobservable errors, $\eps_i$:
% F_n(z)
% = \frac{1}{n} \sum_{i=1}^n \1 (\eps_i \le z),
% \qquad z \in\reals.
%The empirical process corresponding to $F_n$ is
% \alpha_n(z)
% = \sqrt{n} \bigl\{ F_n(z) - F(z) \bigr\},
% \qquad z \in\reals.
%By classical empirical process theory, the latter process converges
%weakly in $\ell^\infty(\reals)$ equipped with the supremum distance:
% \alpha_n \weak\alpha_{F}
% = \BB\circ F,
% \qquad n \to\infty,
%where $\BB$ is a standard Brownian bridge on $[0, 1]$.

%If $F$ is not absolutely continuous, the sequence $\sup_{z \in\reals}
%| \FF_n(z) |$ will no longer be tight, and thus even weak
%hypi-convergence will fail. \js{Conjecture\ldots}

The main arguments underlying the derivation of the limit of $\mathbb
{F}_n$ are as follows.
Let $\mathbb{P}_n=n^{-1} \sum_{i=1}^n \delta_{\mathbf X_i, \varepsilon
_i}$ denote the empirical measure of the sample $(\mathbf{X}_i,
\varepsilon_i)$, $i \in\{1, \ldots, n\}$. For $(z, \bolds{\delta})
\in\bar{\mathbb{R}} \times\mathbb{R}^p$, consider the function
%e6.6 #&#
%
\begin{eqnarray}\label{eqfzdelta}
f_{z, \bolds{\delta}} \dvtx \mathbb{R}^p \times\mathbb{R}& \to& \mathbb{R},
\nonumber\\[-8pt]\\[-8pt]
(\mathbf{x}, \varepsilon) & \mapsto& \mathbh{1} \bigl( \varepsilon\le z + \mathbf
{x}' \bolds{\delta} \bigr),\nonumber
\end{eqnarray}
and let $\mathcal{F}$ denote the collection of all those functions,
that is,
%e6.7 #&#
%
\begin{equation}
\label{eqFclass} \mathcal{F} = \bigl\{ f_{z, \bolds{\delta}}\dvtx  z \in\bar
{\mathbb{R}},
\bolds{\delta} \in\mathbb{R}^p \bigr\}.
\end{equation}
%
% Compactification in the argument $z$ will come in handy later on. For
%$z = -\infty$, we have $f_{z, \vect{\delta}} \equiv0$ and for $z = +
Combining (\ref{eqhateps}) and (\ref{eqhatF}) on the one hand with
(\ref{eqfzdelta}) on the other hand, we find
\[
\hat{F}_n(z) = \frac{1}{n} \sum_{i=1}^n
\mathbh{1} \bigl( \varepsilon_i \le z + \mathbf{X}_i'
( \hat{\bolds{\beta}}_n - \bolds{\beta} ) \bigr) %\\
% &=& \frac{1}{n} \sum_{i=1}^n
% f_{z, \hat{\vect{\beta}}_n - \vect{\beta}}( \vect{X}_i, \eps_i )
=
\mathbb{P}_n f_{z, \hat{\bolds{\beta}}_n - \bolds{\beta}},
\]
where we use the usual operator notation $Qh=\int h \,dQ$ for a signed
measure~$Q$ and a measurable function $h$.
Moreover, let $P$ denote the common law of the random vectors $(\mathbf
{X}_i, \varepsilon_i)$, yielding
$F(z) = \mathbb{E}[ f_{z, \mathbf{0}}( \mathbf{X}_i, \varepsilon_i ) ] = P
f_{z, \mathbf{0}}$ for $z \in\bar{\mathbb{R}}$.
% F(z) = \expec[ f_{z, \vect{0}}( \vect{X}_i, \eps_i ) ] = P f_{z,
Then the empirical process $\mathbb{F}_n$ in (\ref{eqFn}) admits the
decomposition
%e6.8 #&#
%
\begin{eqnarray}
\label{eqFndecomp}
\nonumber
\mathbb{F}_n(z) &=& \sqrt{n} (
\mathbb{P}_n f_{z, \hat{\bolds{\beta}}_n - \bolds{\beta}} - P f_{z,
\mathbf{0}})
\nonumber\\[-8pt]\\[-8pt]
&=& \mathbb{G}_n f_{z, \hat{\bolds{\beta}}_n - \bolds{\beta}} + \sqrt{n} ( P
f_{z, \hat{\bolds{\beta}}_n - \bolds{\beta}} -
P f_{z, \mathbf{0}} ),\nonumber
\end{eqnarray}
where $\mathbb{G}_n$ is shorthand for
$
\mathbb{G}_n = \sqrt{n} (\mathbb{P}_n - P)$.
The decomposition in (\ref{eqFndecomp}) is akin to the one in (\ref
{eqmd}) for the empirical copula process.
If $\hat{\bolds{\beta}}_n$ is consistent for $\bolds{\beta}$, the first
term can be shown to be
\[
\mathbb{G}_n f_{z, \hat{\bolds{\beta}}_n - \bolds{\beta}} = \mathbb{G}_n
f_{z, \mathbf{0}} + o_p(1) = \frac{1}{\sqrt{n}}\sum
_{i=1}^n \bigl\{\mathbh{1}(\varepsilon_i
\le z) - F(z)\bigr\} + o_p(1) % = \alpha_n(z) + o_p(1),
\]
uniformly in $z \in\bar{\mathbb{R}}$. The process on the right-hand
side is the usual empirical process corresponding to $\varepsilon_1,
\ldots, \varepsilon_n$, and its weak convergence is one of the
classical results of empirical process theory.
%The process on the right-hand side is as if the errors $\eps_i$
%themselves were observed, i.e., $\vect{\beta}$ is known, and its weak
%convergence can be handled by standard methods from empirical process
%theory. \ab{References.?}

The treatment of the second term in (\ref{eqFndecomp}) will be based
on a linear expansion of the map $\bolds{\delta} \mapsto P f_{z, \bolds
{\delta}}$ around $\mathbf{0}$. For $(z, \bolds{\delta}) \in\bar{\mathbb
{R}} \times\mathbb{R}^p$, we have
\[
P f_{z, \bolds{\delta}} - P f_{z, \mathbf{0}} = \int_{\mathbb{R}^p} \bigl
\{ F\bigl(z + \mathbf{x}' \bolds{\delta}\bigr) - F(z) \bigr\}
P^X( {d}\mathbf{x} ).
\]
Therefore, if $F$ is continuously differentiable with derivative $f$,
we can expect that
%e6.9 #&#
%
\begin{eqnarray}\label{eqexpansion}
\sqrt{n} \{ P f_{z, \hat{\bolds{\beta}}_n - \bolds{\beta}} - P f_{z,
\mathbf{0}} \} &=& \sqrt{n}
\int_{\mathbb{R}^p} f(z) \mathbf{x}' ( \hat{\bolds{
\beta}}_n - \bolds{\beta} ) P^X ({d}\mathbf{x}) +
o_p(1)
\nonumber\\[-8pt]\\[-8pt]
&=& f(z) \mathbb{E}[\mathbf{X}]' \sqrt{n} ( \hat{\bolds{
\beta}}_n - \bolds{\beta} ) + o_p(1),\nonumber
\end{eqnarray}
which will converge weakly provided $\sqrt n ( \hat{\bolds{\beta}}_n -
\bolds{\beta} )$ converges weakly.
However, if $F$ is not differentiable at a point $z$ or if $f$ exists
but is not continuous in~$z$, then (\ref{eqexpansion}) and as a
consequence weak convergence with respect to the supremum distance may fail.
Still, weak hypi-convergence continues to hold, as the main result of
this section shows.
%Recall our convention of omitting the brackets $[ \cdot]$ when
%working in $\Loc(\TT)$.

%th6.1 #&#
%
\begin{theorem}
\label{theoresweak}
Consider a model of the form (\ref{eqregmod}) such that $(\mathbf{X}_i,
\varepsilon_i)$, $i\in\mathbb{N}$, are i.i.d. random vectors in
$\mathbb{R}^p \times\mathbb{R}$ and such that $\mathbf{X}_i$ and
$\varepsilon_i$ are independent.
Additionally, suppose that the following conditions hold:
\begin{longlist}[(R1)]
\item[(R1)]
The estimator $\hat{\bolds{\beta}}_n$ admits a linear expansion of the form
\[
\sqrt{n} ( \hat{\bolds{\beta}}_n - \bolds{\beta} ) = (
\mathbb{G}_n \psi_1, \ldots, \mathbb{G}_n
\psi_p)' %+ o_p(1)
% = \GG_n \vect{\psi}
+ o_p(1),
\qquad n \to\infty,
\]
in terms of zero-mean, square-integrable functions $\bolds{\psi}_j\dvtx
\mathbb{R}^p \times\mathbb{R}\to\mathbb{R}$, for $j \in\{1,\ldots, p\}$.
\item[(R2)]
The distribution $F$ is absolutely continuous. There exists a version
of its density $f$ which is uniformly bounded and which is \emph{l\`adl\`ag}, that is, which admits right-hand and left-hand limits at
every $z \in\mathbb{R}$:
\[
f(z+) = \lim_{0 < s \to0} f(z + s), \qquad
f(z-) = \lim_{0 < s \to0} f(z - s).
\]
Moreover, $f(\pm\infty):= \lim_{z \to\pm\infty} f(z)=0$.
\item[(R3)]
The norm of $\mathbf X$ is integrable, that is, $\mathbb{E}[ | \mathbf{X} | ]
< \infty$.
% \expec[ | \vect{X} | ] < \infty.
\end{longlist}
Set $\psi=(\psi_1,\ldots, \psi_p)$, $\mathcal{G}= \mathcal{F}\cup
\{ \psi_1, \ldots, \psi_p \} $ and let $\mathbb{G}$ denote a
$P$-Brownian bridge in $\ell^\infty(\mathcal{G})$,
that is, a zero-mean Gaussian process on $\mathcal{G}$ with covariance function
%e6.10 #&#
%
\begin{equation}
\label{eqcovG} \operatorname{cov}( \mathbb{G}g_1,
\mathbb{G}g_2 ) = \operatorname{cov} \bigl( g_1( \mathbf{X},
\varepsilon), g_2( \mathbf{X}, \varepsilon) \bigr), \qquad
g_1, g_2 \in\mathcal{G}.
\end{equation}
Then, in $(L^\infty(\bar\mathbb{R}), {d_{\mathrm{hypi}}} )$, %$(
we have
$ \mathbb{F}_n \rightsquigarrow
\mathbb{F}$ as $n\to\infty$,
% \FF_n \weak
% \FF
% = T \bigl( \GG f_{\cdot, \vect{0}}, \GG\vect{\psi} \bigr)
% = \GG f_{\cdot, \vect{0}} + g( \GG\vect{\psi} )
% , \qquad n \to\infty,
where the limiting process $\mathbb{F}$ can be written as $\mathbb
{F}(\pm\infty) = 0$ a.s. and
%e6.11 #&#
%
\begin{eqnarray}
\label{eqFF} \mathbb{F}(z) &=& \mathbb{G}f_{z, \mathbf{0}}  - f(z-)
\int
_{-\infty}^0 P^X\bigl( \bigl\{ \mathbf{x}\dvtx
\mathbf{x}' \mathbb{G}\bolds{\psi} < y \bigr\} \bigr) \,{d}y
\nonumber\\[-8pt]\\[-8pt]
&&{} + f(z+) \int_0^{+\infty}
P^X\bigl( \bigl\{ \mathbf{x}\dvtx  \mathbf{x}' \mathbb{G}\bolds{\psi}
> y \bigr\} \bigr) \,{d}y, \qquad z \in\mathbb{R}.\nonumber
\end{eqnarray}
\end{theorem}

Note that the limit in (\ref{eqFF}) is not c\`adl\`ag, whence the
classical Skorohod-topologies cannot be applied in the present context.

% Assumption (R1) postulates an asymptotically linear expansion of $
The influence function $\bolds{\psi}=(\bolds\psi_1,\ldots, \bolds\psi_p)$
in (R1) depends on the estimator and on the true model.
% , i.e., on the parameter vector $\vect{\beta}$ and the laws of $
A classical example is given by the ordinary least squares estimator:
if the errors $\varepsilon_i$ have mean zero and finite variance and
if the components of $\mathbf{X}$ have finite second moments and the $p
\times p$ matrix $\mathbb{E}( \mathbf{X} \mathbf{X}' )$ is invertible, then
\begin{eqnarray*}
\sqrt{n} ( \hat{\bolds{\beta}}_n - \bolds{\beta} ) % &=& \biggl\{ \frac{1}{n}
&=&
\frac{1}{\sqrt{n}} \sum_{i=1}^n \bigl\{
\mathbb{E}\bigl( \mathbf{X} \mathbf{X}' \bigr) \bigr\}^{-1}
\mathbf{X}_i \varepsilon_i + o_p(1).
\end{eqnarray*}
%
% General classes of estimators that admit such an expansion include M-
%and Z-estimators, see for instance Chapters~3.2 and~3.3 in
% \ab{Including ordinary least squares as we use this in the example
%below.} \sv{spell out the expansion for least squares?}

If $f$ happens to be continuous in $z$, then $f(z-) = f(z+) = f(z)$,
and we obtain that
$
\mathbb{F}(z) = \mathbb{G}f_{z, \mathbf{0}} + f(z) \mathbb{E}[ \mathbf
{X}' ] \mathbb{G}\bolds{\psi}$,
which, under (R1), coincides with the limit of the classical
representation in (\ref{eqreprffn}). If $f$ is continuous everywhere,
then $\mathbb{F}$ is almost surely continuous, and, by Corollary~\ref
{corK} with $K = \bar{\mathbb{R}}$, the weak convergence in $\mathbb
{F}_n \rightsquigarrow\mathbb{F}$ takes place with respect to the
supremum distance.

%continuous function only yields locally uniform convergence, i.e., in $
%replace $\reals$ by $[-\infty, +\infty]$; for $z \in\{-\infty, +\infty
%additional thought.}

%ex6.2 #&#
%
\begin{example}[(Mixtures of exponential distributions)]%
Consider the probability density function
%e6.12 #&#
%
\begin{equation}
\label{eqdoubleExponential} f_{\bolds{\theta}}(z) = \cases{ \displaystyle
\frac{1}{\theta_-}
e^{z/\theta_-}, &\quad if $z < 0$,
\vspace*{5pt}\cr
\displaystyle(1-w) \frac{1}{\theta_+}
e^{-z/\theta_+}, &\quad if $z > 0$,}
\end{equation}
%
% f_{\vect{\theta}}(z)
% =
% \cases{
% \frac{\theta_+}{\theta_-} \frac{1}{\theta_+ + \theta_-} e^{z/
% \frac{\theta_-}{\theta_+} \frac{1}{\theta_+ + \theta_-} e^{-z/
% \end{cases}
where $w = \theta_+/(\theta_- + \theta_+)$.
This density is a mixture of the exponential distribution on $(-\infty, 0)$ with mean $-\theta_-$ and the exponential distribution on $(0,
\infty)$ with mean~$\theta_+$, with weights chosen so that the total
mean is zero.
%The mean of the distribution is $-w \theta_- + (1-w) \theta_+$, which
%is equal to zero if
% w = \frac{\theta_+}{\theta_- + \theta_+},
%which we will assume throughout. The corresponding distribution
%function is given by
% F_{\vect{\theta}}(x) = \int_{-\infty}^x f_{\vect{\theta}}(z) \diff
%z =
% \cases{
% w e^{z/\theta_-} & \mbox{if $z \le0$,} \\
% w + (1-w) (1 - e^{-z/\theta_+}) & \mbox{if $z \ge0$}
% \end{cases}
The left-hand and right-hand limits of $f_{\bolds{\theta}}$ at $0$ are
\[
f_{\bolds{\theta}}(0-) = \frac{\theta_+}{\theta_-} \frac
{1}{\theta_- + \theta_+}, \qquad
f_{\bolds{\theta}}(0+) = \frac{\theta_-}{\theta_+} \frac
{1}{\theta_- + \theta_+}.
\]
If $\theta_-$ is different from $\theta_+$, these limits are
different, and thus the associated distribution function, $F_{\bolds
{\theta}}$, is not continuously differentiable at $0$.
See the left-hand side of Figure~\ref{figpdf-cdf} for the graph of
$F_\theta$ when $(\theta_-, \theta_+)$ is equal to $(1,
4)$.\looseness=-1
%See the upper row of Figure~\ref{figpdf-cdf} for the graphs of $F_

%f2 #&#
%
\begin{figure}%[ht!]

\includegraphics{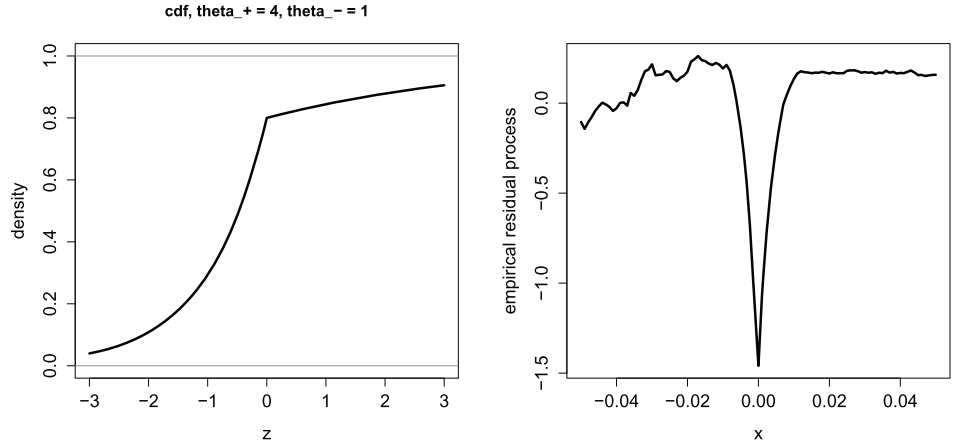}

\caption{\textup{Left}: cumulative distribution function of the mixed
double exponential distribution in~(\protect\ref{eqdoubleExponential}).
\textup{Right}: trajectories of the corresponding
empirical residual process $\mathbb{F}_n$ for $n=10^6$. In both cases:
$(\theta_-, \theta_+)=(1,4)$.%
% Left: $\theta_+ = \Sexpr{betaposA}$, $\theta_- = \Sexpr{betanegA}$.
%Right: $\theta_+ = \Sexpr{betaposB}$, $\theta_- = \Sexpr{betanegB}$.}
}\label{figpdf-cdf}
\end{figure}

Now, consider the linear regression model in (\ref{eqregmod}) with
$p=1$, $X_i\sim N(0,1)$ independent of $\varepsilon_i$, and with
$\varepsilon_i$ distributed according to (\ref
{eqdoubleExponential}). %with \eqref{eqmeanZero}.
The parameter $\bolds{\beta} \in\mathbb{R}$ is estimated by
ordinary
least squares, $\hat{\bolds{\beta}}_n$, and the corresponding empirical\vadjust{\goodbreak}
residual process $\mathbb{F}_n$ is calculated as in (\ref{eqFn}).
Theorem~\ref{theoresweak} implies that $\mathbb{F}_n$ converges in
$(L^\infty( \bar\mathbb{R}), {d_{\mathrm{hypi}}} )$ to the
process $\mathbb{F}$ given by (note the simplification arising from
$\mathbb{E}[X] = 0$)
\[
\mathbb{F}(z) = \mathbb{G}f_{z, \mathbf{0}} + \bigl( f_{\bolds{\theta
}}(z+) -
f_{\bolds{\theta}}(z-) \bigr) \int_0^{+\infty}
P^X\bigl( \{ x\dvtx  x \mathbb{G}\psi> y \} \bigr) \,{d}y,
\]
where $\mathbb{G}$ is a $P$-Brownian bridge for $P=P^X \otimes
P^\varepsilon$ and where $f_{z, \mathbf{0}}$ and $\psi$ are certain
functions in $L^2(P)$. We find that $\mathbb{F}(z) = \mathbb{G}f_{z,
\mathbf{0}}$ for $z \ne0$, a continuous Gaussian process. The only
discontinuity occurs at $z = 0$, when the left-hand and right-hand
limits of $f_{\bolds{\theta}}$ are different. The ``spike'' in $\mathbb
{F}_n$ then goes upward or downward according to whether $f_{\bolds
{\theta}}(z+)$ is larger than or smaller than $f_{\bolds{\theta}}(z-)$.
A simulated typical trajectories of $\mathbb{F}_n(z)$ for $n=10^6$ and
$z\in[-0.05, 0.05]$
is shown on the right-hand side of Figure~\ref{figpdf-cdf} when
$(\theta_-, \theta_+)$ is equal to $(1, 4)$.
\end{example}

\begin{appendix}
% \section{Auxiliary results}\label{appaux}
%
% This appendix contains a number of auxiliary results, the proofs of
%which are all deferred to the supplement [
%some of the Lemmas completely to the supplement? Keep only those which
%are cited throughout the main file?} \js{Give a brief outline of the
%subsections of this appendix. Or even split the appendix in three
%appendices, one for each subsection?}
%s7 #&#
\section{Verifying hypi-convergence}\label{APPVERIF}

In this appendix, we provide some tools for showing convergence of a
sequence of functions with respect to the hypi-semimetric. Proofs are
deferred to Appendix~D in the supplement
[\citet{BucSegVol14supp}].

%s7.1 #&#
\subsection{Pointwise convergence and convergence of sums}\label{appsum}

Let $(\mathbb{T}, d)$ be a metric space. For $f\dvtx  \mathbb{T}\to
\mathbb{R}$, define extended real-valued functions ${f_\wedge} $ and
${f_\vee} $ as in~(\ref{eqfinf}) and (\ref{eqfsup}), respectively.
Since we do not require $f$ to be locally bounded, ${f_\wedge} $ and
${f_\vee} $ can attain $-\infty$ and $+\infty$, respectively.

For $f_n\dvtx  \mathbb{T}\to\mathbb{R}$, we say that $f_n$ epi-converges
to $\alpha\in\mathbb{R}$ at $x \in\mathbb{T}$ if the following two
conditions are met:
%e7.1 #&#
%
\begin{equation}
\label{eqdefepi} \cases{\displaystyle\mbox{\phantom{i}(i) } \forall x_n \to x\dvtx  \liminf
_{n \to\infty} f_n (x_n) \ge\alpha,
\vspace*{5pt}\cr
\displaystyle\mbox{(ii) }\exists x_n \to x\dvtx  \limsup_{n \to\infty}
f_n(x_n) \le\alpha.}
\end{equation}
Similarly, $f_n$ hypo-converges to $\alpha$ at $x$ if
%e7.2 #&#
%
\begin{equation}
\label{eqdefhypo} \cases{ \displaystyle\mbox{\phantom{i}(i) } \forall x_n \to x\dvtx  \limsup
_{n \to\infty} f_n (x_n) \le\alpha,
\vspace*{5pt}\cr
\displaystyle\mbox{(ii) } \exists x_n \to x\dvtx  \liminf_{n \to\infty}
f_n(x_n) \ge\alpha,}
\end{equation}
which is equivalent to epi-convergence of $-f_n$ to $-\alpha$ at $x$.
If additionally \mbox{$f\dvtx \mathbb{T}\to\mathbb{R}$}, then $f_n$ is said to
epi- or hypo-converge to $f$ at $x$ if $\alpha=f(x)$ in the preceding
conditions. According to Proposition~\ref{prophypiconv}, $f_n$
hypi-converges to $f$ in $\ell_{\mathrm{loc}}^\infty(\mathbb{T})$
if and only if $f_n$ epi-converges to ${f_\wedge} $ and hypo-converges
to ${f_\vee} $ at every $x \in\mathbb{T}$. For $x \in\mathbb{T}$
and $\varepsilon> 0$, let $B(x, \varepsilon) = \{ y \in\mathbb{T}\dvtx
d(x, y) < \varepsilon\}$.

\begin{lemma}[(Convergence of hulls)]\label{lemconthull}
Let $f_n\dvtx \mathbb{T}\to\mathbb{R}$, $x\in\mathbb{T}$ and $\alpha
\in\mathbb{R}$. Then
$f_n$ epi-converges to $\alpha$ at $x$ if and only if $f_{n, \wedge}$
epi-converges to $\alpha$ at $x$, and $f_n$ hypo-converges to $\alpha
$ at $x$ if and\vadjust{\goodbreak} only if $f_{n, \vee}$ hypo-converges to $\alpha$ at~$x$. Moreover,
%e7.3 #&#
%
\begin{equation}
\label{eqhull1} f_n(x_n) \to\alpha\qquad\forall
x_n \to x
\end{equation}
is equivalent to
%e7.4 #&#
%
\begin{equation}
\label{eqhull2} f_{n,\wedge}(x_n) \to\alpha\quad\mbox{and}\quad
f_{n,\vee}(x_n) \to\alpha\qquad\forall x_n \to x.
\end{equation}
\end{lemma}

The following three lemmas contain results on hulls of sums and on
epi-, hypo- and hypi-convergence of sums.

%le7.2 #&#
%
\begin{lemma}[(On sums of hulls and hulls of sums)]\label{lemsumsandwich}
For $f, g\dvtx  \mathbb{T}\to\mathbb{R}$ such that $g_\wedge$ and
$g_\vee$ are both finite, we have
\begin{eqnarray*}
{f_\wedge} + {g_\wedge} &\le& (f+g)_\wedge
\le{f_\wedge} + {g_\vee},
\\
{f_\vee} + {g_\wedge} &\le& (f+g)_\vee
\le{f_\vee} + {g_\vee}.
\end{eqnarray*}
In particular, if $g$ is continuous in $x \in\mathbb{T}$, then
$(f+g)_{\wedge}(x) = {f_\wedge} (x) + g(x)$ and $(f+g)_{\vee}(x) =
{f_\vee} (x) + g(x)$.
\end{lemma}

%le7.3 #&#
%
\begin{lemma}[(Epi- and hypo-convergence of hulls of sums)]\label{lemsumcont}
Let $f_n,g_n\dvtx \mathbb{T}\to\mathbb{R}$ and let $x \in\mathbb{T}$ be
such that $g_n(x_n) \to\beta\in\mathbb{R}$ for all sequences $x_n
\to x$. If $f_{n,\wedge}$ epi-converges to $\alpha$ at $x$, then
$(f_n+g_n)_\wedge$ epi-converges to $\alpha+\beta$ at $x$. Similarly
for upper semicontinuous hulls and hypo-convergence.
\end{lemma}

%le7.4 #&#
%
\begin{lemma}[(Hypi-convergence of sums)]\label{lemsum}
Let $\mathbb{T}$ be locally compact and separable. If $f_n$ and $g_n$
hypi-converge to $f$ and $g$ in $\ell_{\mathrm{loc}}^\infty(\mathbb
{T})$, respectively, and if at every point $x \in\mathbb{T}$, at least
one of the two functions $f$ or $g$ is continuous, then $f_n+g_n$
hypi-converges to $f+g$.
\end{lemma}

%s7.2 #&#
\subsection{Upper and lower semicontinuous extensions}\label{appext}

The limit processes in Theorems~\ref{theocopweak} and~\ref
{theotailweak} are defined by extending a continuous function defined
on a dense subset of a metric space to the whole space. In this
section, some useful elementary properties of such extensions are
recorded. The main tool is Corollary~\ref{corextension}, giving a
criterion for proving hypi-converg\-ence to a function defined by such
an extension procedure.

Let $(\mathbb{T}, d)$ be a metric space, let $A \subset\mathbb{T}$
be dense, and let $f\dvtx  A \to\mathbb{R}$. Extend the domain of $f$
from $A$ to the whole of $\mathbb{T}$ by
%e7.5 #&#
%e7.6 #&#
%
\begin{eqnarray}
\label{eqfinfA} {f_\wedge^{ A\dvtx \mathbb{T} } } (x) &=& \sup
_{\varepsilon> 0} \inf f\bigl( B(x, \varepsilon) \cap A \bigr) \in
[-\infty,
\infty],
\\
\label{eqfsupA} {f_\vee^{ A\dvtx \mathbb{T} } } (x) &=& \inf
_{\varepsilon> 0} \sup f\bigl( B(x, \varepsilon) \cap A \bigr) \in
[-\infty,
\infty],
\end{eqnarray}
for $x \in\mathbb{T}$; as before, $B(x, \varepsilon) = \{ y \in
\mathbb{T}\dvtx  d(x, y) < \varepsilon\}$ is the open ball centered at $x$
of radius $\varepsilon$.
Note that these definitions also make sense if $A = \mathbb{T}$, and
that for $f \in\ell_{\mathrm{loc}}^\infty(\mathbb{T})$ we have
${f_\wedge^{ \mathbb{T}\dvtx \mathbb{T} } } = {f_\wedge} $ and ${f_\vee
^{ \mathbb{T}\dvtx \mathbb{T} } } = {f_\vee} $; see the definitions in
(\ref{eqfinf})~and~(\ref{eqfsup}).\vadjust{\goodbreak}

Clearly, ${f_\wedge^{ A\dvtx \mathbb{T} } } (x) \le f(x) \le{f_\vee^{
A\dvtx \mathbb{T} } } (x)$ for every $x \in A$.
For any open set $U \subset\mathbb{T}$, we have
\[
\inf{f_\wedge^{ A\dvtx \mathbb{T} } } (U) = \inf f(U \cap A), \qquad
\sup{f_\vee^{ A\dvtx \mathbb{T} } } (U) = \sup f(U \cap A).
\]
The functions ${f_\wedge^{ A\dvtx \mathbb{T} } } $ and ${f_\vee^{
A\dvtx \mathbb{T} } } $ from $\mathbb{T}$ into $[-\infty,+\infty]$ are
lower and upper semicontinuous, respectively. If every $x$ in $A$
admits a neighborhood on which~$f$ is bounded, then ${f_\wedge^{
A\dvtx \mathbb{T} } } $ and ${f_\vee^{ A\dvtx \mathbb{T} } } $ are real-valued.

If\vspace*{1pt} $f$ is continuous at $x \in A$, then ${f_\wedge^{ A\dvtx \mathbb{T} } }
(x) = {f_\vee^{ A\dvtx \mathbb{T} } } (x) = f(x)$, and ${f_\wedge^{
A\dvtx \mathbb{T} } } $ and ${f_\vee^{ A\dvtx \mathbb{T} } } $, seen as
functions on $\mathbb{T}$, are continuous at $x$, too.
The following lemma shows that, if $f$ is continuous on the whole of
$A$, then its domain does not really matter insofar as the extension is
concerned.
%It suggest moreover that Belgium is going to win the 2014 World Cup in
%Brazil.

%le7.5 #&#
%
\begin{lemma} \label{lemextext}
Let $E \subset A \subset\mathbb{T}$ be such that $E$ is dense in
$\mathbb{T}$. Let $f\dvtx  A \to\mathbb{R}$ and consider the restriction\vspace*{1pt}
$f|_E\dvtx  E \to\mathbb{R}$ of $f$ to $E$ and the extensions
$(f|_E)_\wedge^{E\dvtx \mathbb{T}}$ and $(f|_E)_\vee^{E\dvtx \mathbb{T}}$ of
$f|_E$ to $\mathbb{T}$. If $f$ is continuous,\vspace*{1pt} then $(f|_E)_\wedge
^{E\dvtx \mathbb{T}} = {f_\wedge^{ A\dvtx \mathbb{T} } } $ and $(f|_E)_\vee
^{E\dvtx \mathbb{T}} = {f_\vee^{ A\dvtx \mathbb{T} } } $.
\end{lemma}

The following two results provide criterions for proving epi-, hypo- or
hypi-converg\-ence to a semicontinuous extension.

%pr7.6 #&#
%
\begin{proposition}
\label{propextension}
Let $A \subset\mathbb{T}$ be dense and let $f\dvtx  A \to\mathbb{R}$ be
continuous.
%and let $\finf$ be its lower semi-continuous extension from $A$ to $
Assume that ${f_\wedge^{ A\dvtx \mathbb{T} } } $ is real-valued.\vspace*{1.5pt} If the
functions $f_n\dvtx  \mathbb{T}\to\mathbb{R}$ converge pointwise on $A$
to $f$ and if $\liminf_n f_n(x_n) \ge{f_\wedge^{ A\dvtx \mathbb{T} } }
(x)$ whenever\vspace*{1.5pt} $x_n \in\mathbb{T}$ converges to $x \in\mathbb{T}$,
then $f_n$ epi-converges to ${f_\wedge^{ A\dvtx \mathbb{T} } } $.
Similarly for hypo-convergence to~${f_\vee^{ A\dvtx \mathbb{T} } } $.
\end{proposition}

%co7.7 #&#
%
\begin{corollary}
\label{corextension}
Let $A \subset\mathbb{T}$ be dense. Let $f\dvtx  A \to\mathbb{R}$ be
continuous and suppose that its lower and upper semicontinuous
extensions %from $A$ to~$\TT$
${f_\wedge^{ A\dvtx \mathbb{T} } } $ and ${f_\vee^{ A\dvtx \mathbb{T} } } $
are real-valued. Let $f^*\dvtx  \mathbb{T}\to\mathbb{R}$ be such that
${f_\wedge^{ A\dvtx \mathbb{T} } } \le f^* \le{f_\vee^{ A\dvtx \mathbb{T} }
} $. Then ${f_\wedge^{ A\dvtx \mathbb{T} } } = (f^*)_\wedge$ and ${f_\vee
^{ A\dvtx \mathbb{T} } } = (f^*)_\vee$ on $\mathbb{T}$. Moreover, if the
functions $f_n\dvtx  \mathbb{T}\to\mathbb{R}$ are locally bounded and verify
\[
\forall x \in\mathbb{T}\dvtx  \forall x_n \to x\dvtx  {f_\wedge^{ A\dvtx \mathbb{T} }
} (x) \le\liminf_{n \to\infty} f_n(x_n) \le
\limsup_{n \to\infty} f_n(x_n) \le
{f_\vee^{ A\dvtx \mathbb{T} } } (x),
\]
then $f_n$ hypi-converges to $f^*$.
\end{corollary}

%s8 #&#
\section{Weak convergence and semimetric spaces}\label{APPWEAKSEMI}

The workhorses of the theory of weak convergence in metric spaces are
the continuous mapping theorem, the extended continuous mapping
theorem, and for normed vector spaces, the functional delta method;
see, for instance, Theorems~1.3.6, 1.11.1 and 3.9.4 in \citet
{vandwell1996}. Thanks to these theorems, weak convergence of many
empirical processes can be shown and\vadjust{\goodbreak} can be exploited to conclude weak
convergence of sequences of appropriately normalized estimators and
test statistics. However, the space of interest in this paper, $(\ell
_{\mathrm{loc}}^\infty(\mathbb{T}), {d_{\mathrm{hypi}}} )$; see
Proposition~\ref{prophypiconv}, is not a metric space but rather a
semimetric space. Moreover, addition of functions is ill-compatible
with the hypi-semimetric: if $f \in\ell_{\mathrm{loc}}^\infty
(\mathbb{T})$ is not continuous, then ${d_{\mathrm{hypi}}} (f + g,
0)$ need not be equal to zero even if ${d_{\mathrm{hypi}}} (g, -f) =
0$. Hence, addition is not well defined on the space of equivalence
classes of functions at hypi-distance zero.

In this appendix, versions of the (extended) continuous mapping theorem
and the functional delta method are given that are adapted to
semimetric spaces. In particular, the maps under consideration are not
required to be defined on equivalence classes of points at distance
zero but rather on the original semimetric space itself.
% This kind of generality is needed to be able to extend Slutsky's
%lemma to weak hypi-convergence (Lemma~\ref{lemSlutsky}).
Proofs are deferred to Appendix~E in the
supplement [\citet{BucSegVol14supp}].

% The continuous mapping theorem is one of the workhorses in weak
%convergence theory [\citet[][Theorem~1.3.6]{vandwell1996}. We
%formulate and prove a version that is adapted to semimetric spaces
%(Theorem~\ref{thmcmt}). In particular, the map under consideration is
%not required to be defined on equivalence classes of points at
%distance zero but rather on the original semimetric space itself. This
%allows to extend Slutsky's lemma to hypi-convergence (Lemma~
%theorem and the delta method. Split into two or more subsections as
%well?}

Let $(\mathbb{D},d)$ be a semimetric space. For $x \in\mathbb{D}$,
put $[x] = \operatorname{cl}\{x\}$, the set of $y \in\mathbb{D}$
such that $d(x,y) = 0$. Since $d(x', y') = d(x, y)$ whenever $x' \in
[x]$ and $y' \in[y]$, we can, abusing notation, define a metric
$d([x],[y]):= d(x,y)$ on the quotient space $[\mathbb{D}] = \{ [x]\dvtx
x \in\mathbb{D}\}$. %(Abuse of notation\dvtx  use `$d$' for semimetric on $
Let $[ \cdot]$ denote the map $\mathbb{D}\to[\mathbb{D}]\dvtx  x
\mapsto[x]$. Obviously, $[ \cdot]$ is continuous. The image of an
open (closed) subset of $\mathbb{D}$ under $[ \cdot]$ is open
(closed) in $[\mathbb{D}]$.

% Let $\mathcal{D}$ be the collection of subsets $A$ of $\DD$ with the
%property that $x \in A$ implies $[x] \subset A$. The collection $
%complementation and arbitrary unions. Since $\mathcal{D}$ contains the
%closed sets, it also contains the open sets and even the Borel $

Let $\mathcal{B}(\mathbb{D})$ and $\mathcal{B}([\mathbb{D}])$ be
the Borel $\sigma$-fields on $(\mathbb{D},d)$ and $([\mathbb
{D}],d)$, respectively, that is, the smallest $\sigma$-fields
containing the open sets. There is a one-to-one correspondence between
both $\sigma$-fields:  for $B \in\mathcal{B}(\mathbb{D})$, the set
$[B] = \{ [x]\dvtx  x \in B \}$ is a Borel set in $\mathbb{D}$, and
conversely, every Borel set $B$ of $\mathbb{D}$ can be written as
$\bigcup_{[x] \in[B]} [x]$; in particular $x\in B$ if and only if
$[x] \in[B]$. A Borel law $L$ on $(\mathbb{D}, \mathcal{B}(\mathbb
{D}))$ induces a Borel law $L \circ[ \cdot]^{-1}$ on $([\mathbb
{D}], \mathcal{B}([\mathbb{D}]))$ and vice versa.

One of the merits of Hoffman--J{\o}rgensen weak convergence is that
measurability requirements are relaxed. In the context of semimetric
spaces, measurability issues require, perhaps, some extra care.
% This is the purpose of the following lemma.

%le8.1 #&#
%
\begin{lemma}[(Measurability)]\label{lemmeasurable}
Let $(\mathbb{D},d)$ and $(\mathbb{E},e)$ be semimetric spaces. Let
$g\dvtx  \mathbb{D}\to\mathbb{E}$ be arbitrary. Then the set $D_g$ of $x
\in\mathbb{D}$ such that $g$ is not continuous in $x$ is Borel
measurable. More generally, $g^{-1}(B) \setminus D_g$ is a Borel set in
$\mathbb{D}$ for every Borel set $B$ in $\mathbb{E}$.
\end{lemma}

In our version of the continuous mapping theorem, the map $g$ is
defined on~$\mathbb{D}$ and not on $[\mathbb{D}]$, that is, even if
$d(x, y) = 0$, it may occur that $g(x) \ne g(y)$. Therefore, we cannot
directly apply Theorem~1.3.6 in \citet{vandwell1996}.
Nevertheless, the proof is inspired from the proof of that theorem.

%th8.2 #&#
%
\begin{theorem}[(Continuous mapping)]
\label{thmcmt}
Let $(\mathbb{D},d)$ be a semimetric space and let $(\mathbb{E},e)$
be a metric space. Let $g\dvtx  \mathbb{D}\to\mathbb{E}$ be arbitrary
and let $D_g$ be the set of $x \in\mathbb{D}$ such that $g$ is not
continuous in $x$. Let $(\Omega_\alpha, \mathcal{A}_\alpha,
P_\alpha)$, $\alpha\in A$, be a net of probability spaces and let
$X_\alpha\dvtx  \Omega_\alpha\to\mathbb{D}$ be arbitrary maps; let
$(\Omega, \mathcal{A}, P)$ be a probability space and let $X\dvtx  \Omega
\to\mathbb{D}$ be Borel measurable. If $[X_\alpha] \rightsquigarrow
[X]$ in $[\mathbb{D}]$ and if $X(\Omega) \subset\mathbb{D}\setminus
D_g$, then $g(X)\dvtx  \Omega\to\mathbb{E}$ is Borel measurable and
$g(X_\alpha) \rightsquigarrow g(X)$ in $\mathbb{E}$.
\end{theorem}

%
%The continuous mapping theorem can be used to derive a version of
%Slutsky's Lemma for weak hypi-convergence in $\Loc(\TT)$. \ab{I wonder
%whether we should state this Lemma in the main section on weak
%hypi-convergence?}
%
%Let $(\TT,\dT)$ be a locally compact, separable metric space. Let
%$X_n, Y_n\dvtx  \Omega_n \mapsto\loc(\TT)$ be arbitrary maps and let $X\dvtx
%hypi-semimetric. If $[X_n] \dto[X]$ and $[Y_n] \dto[0]$ in $\Loc(
%
%Let $\DD$ be the product space $\loc(\TT) \times\loc(\TT)$ equipped
%with the box semimetric $d_{\DD}((f_1,f_2), (f_1',f_2')) = \max\{
%Consider the map $g\dvtx  \DD\to\EE$ given by $g(f_1,f_2) = [f_1 + f_2]$.
%
%By Lemma~\ref{lemsum} and by continuity of the map $\loc(\TT) \mapsto
%continuous at pairs $(f_1, f_2)$ such that $f_2$ is continuous, in
%particular if $f_2 = 0$. In view of the Continuous Mapping Theorem~
%0)]$ in $[\DD]$, where $[(f_1,f_2)]$ is the equivalence class of pairs
%$(f_1',f_2')$ in $\DD$ at distance zero from $(f_1,f_2)$, and where $[
%with the natural metric induced by $d_{\DD}$.
%
%The distance (in $[\DD]$) between $[(X_n, Y_n)]$ and $[(X_n, 0)]$ is
%equal to the distance (in $\DD$) between $(X_n, Y_n)$ and $(X_n, 0)$,
%which is equal to $\dhypi(Y_n, 0) = \dhypi([Y_n], [0])$. Since $[Y_n]
%see Lemma~1.10.2(iii) in \citet{vandwell1996}. By item (i) of the
%same lemma, the desired conclusion $[(X_n, Y_n)] \dto[(X, 0)]$ is a
%consequence from $[(X_n, 0)] \dto[(X, 0)]$, to be proven. But the
%latter convergence follows from the Continuous Mapping Theorem~
%[(f, 0)]$.
%}
%

In many circumstances,
%Often,
one needs a refined version of the continuous mapping theorem that
covers maps $g_n(X_n)$, rather than $g(X_n)$ for a fixed $g$. The
following statement and proof are inspired from Theorem~1.11.1(i) and
Problem~1.11.1 in \citet{vandwell1996} and Theorem~18.11(i) in
\citet{vdvaart1998}.

%th8.3 #&#
%
\begin{theorem}[(Extended continuous mapping)]\label{thmextendedContinuousMapping}
Let $(\mathbb{D}, d)$ be a semimetric space and let $(\mathbb{E}, e)$
be a metric space. For integer $n \ge0$, let there be probability
spaces $\Omega_n$, subsets $\mathbb{D}_n \subset\mathbb{D}$, maps
$X_n\dvtx  \Omega_n \to\mathbb{D}_n$ and maps $g_n\dvtx  \mathbb{D}_n \to
\mathbb{E}$. Assume the following two conditions:
\begin{itemize}
\item
For every $x_0 \in\mathbb{D}_0$ and for every subsequence
$(x_{n_k})_k$ with $x_{n_k} \in\mathbb{D}_{n_k}$ for all $k$ and such
that $x_{n_k} \to x_0$ as $k \to\infty$, we have $g_{n_k}(x_{n_k})
\to g_0(x_0)$.
\item
The map $X_0$ is Borel measurable and $[X_n] \rightsquigarrow[X_0]$ in
$([\mathbb{D}], d)$.
\end{itemize}
If $g_0(X_0)$ is Borel measurable, then $g_n(X_n) \rightsquigarrow
g_0(X_0)$ in $(\mathbb{E}, e)$. If $g_0(X_0)$ is not Borel measurable,
there still exists a version $X_0'$ of $X_0$ such that $g_0(X_0')$ is
Borel measurable, and thus $g_n(X_n) \rightsquigarrow g_0(X_0')$ in
$(\mathbb{E}, e)$.
\end{theorem}

%ad8.4 #&#
%
\begin{addendum}
\label{addextendedContinuousMapping}
The law of $X_0$ is concentrated on the set
\[
\mathbb{D}_\infty= \bigcap_{k \ge1}
\operatorname{cl} \biggl( \bigcup_{m \ge k}
\mathbb{D}_m \biggr) = \limsup_{n \to\infty}
\mathbb{D}_n,
\]
which is closed in $\mathbb{D}$. The restriction of the map $g_0$ to
$\mathbb{D}_0 \cap\mathbb{D}_\infty$ is continuous. Whether
$g_0(X_0)$ is measurable or not, there always exists a version $X_0'$
of $X_0$ which takes values in $\mathbb{D}_0 \cap\mathbb{D}_\infty$
and for which $g_0(X_0')$ is Borel measurable.
\end{addendum}

%co8.5 #&#
%
\begin{corollary}
\label{corextendedContinuousMapping}
If in Theorem~\ref{thmextendedContinuousMapping}, $(\mathbb{E}, e)$
is a semimetric space rather than a metric space, the conclusion still
holds with $[g_n(X_n)]$ converging weakly to $[g_0(X_0)]$ or
$[g_0(X_0')]$, respectively, in $([\mathbb{E}], e)$.
\end{corollary}

% The proofs are given in the supplementary file
% Note that,
% for metric rather than semimetric spaces, there exist several
%versions of the extended continuous mapping theorem\dvtx  see
%Theorem~1.11.1 and Problem~1.11.1 in \citet{vandwell1996} and
%Theorem~18.11 in \citet{vdvaart1998}. The variants differ in the
%following points\dvtx
% \begin{itemize}
% \item
% The law of $X_0$ may or may not be assumed to be separable.
% \item
% The map $g_0(X_0)$ may or may not be assumed to be Borel measurable.
% \item
% The convergence condition on $g_n$ may be relaxed to sequences rather
%than subsequences, i.e., it is required that whenever $x_n \in\DD_n$
%for all integer $n \ge1$ and $x_n \to x \in\DD_0$, we have $g_n(x_n)
% \end{itemize}
% Theorem~\ref{thmextendedContinuousMapping} combines these variants
%in a way which will be suitable for the derivation of a functional
%delta method in semimetric vector spaces.

The formulation of Theorem~\ref{thmextendedContinuousMapping} has
been chosen to make it suitable for establishing a variant of the
functional delta method in semimetric vector spaces. In the remaining
part of this section, let $\mathbb{D}$ and $\mathbb{E}$ be real
vector spaces and let $d$ and $e$ be semimetrics on $\mathbb{D}$ and
$\mathbb{E}$, respectively. Addition is not required to be continuous.
Worse still, addition need not even be well defined on equivalence classes.

%de8.6 #&#
%
\begin{defin}[(Semi-Hadamard differentiability)]
Let $\Psi\dvtx  \mathbb{D}_\Psi\to\mathbb{E}$ with $\mathbb{D}_\Psi
\subset\mathbb{D}$. Let $x \in\mathbb{D}_\Psi$ and $\mathbb
{W}\subset\mathbb{D}$. Then $\Psi$ is said to be semi-Hadamard
differentiable at $x$ tangentially to $\mathbb{W}$ if there exists a
map $d \Psi_x\dvtx  \mathbb{W}\to\mathbb{E}$, called the
(semi-)derivative of $\Psi$ at $x$, with the following property: for
every $w \in\mathbb{W}$, every sequence $( t_n )_n$ in $(0, \infty)$
such that $t_n \to0$ and every sequence $( w_n )_n$ in $\mathbb{D}$
such that $x + t_n w_n \in\mathbb{D}_\Psi$ for all $n$ and $w_n \to
w$ as $n \to\infty$, we have
\[
t_n^{-1} \bigl( \Psi(x + t_n w_n)
- \Psi(x) \bigr) \to d\Psi_x(w), \qquad n \to\infty.
\]
\end{defin}

The derivative $d\Psi_x$ is not assumed to be linear or continuous.
Still, in Addendum~\ref{adddelta} below, we will see that $d\Psi_x$
does enjoy some kind of continuity property. An extension of the chain
rule similar to Lemma~3.9.3 in \citet{vandwell1996} is
straightforward and is therefore omitted.

%th8.7 #&#
%
\begin{theorem}
\label{thmdelta}
Let $\mathbb{D}_\Psi\subset\mathbb{D}$ and let $\Psi\dvtx  \mathbb
{D}_\Psi\to\mathbb{E}$ be semi-Hadamard differentiable at $x \in
\mathbb{D}_\Psi$ tangentially to $\mathbb{W}\subset\mathbb{D}$
with derivative $d \Psi_x$. Let $Y_n$, $n \ge1$, and $X$ be maps from
probability spaces into $\mathbb{D}$ such that $Y_n$ takes values in
$\mathbb{D}_\Psi$ and $X$ takes values in $\mathbb{W}$. Assume that
$X$ is Borel measurable and that there exists a positive sequence $r_n$
tending to infinity such that, in $([\mathbb{D}], d)$,
\[
\bigl[r_n (Y_n - x)\bigr] \rightsquigarrow[X], \qquad n
\to\infty.
\]
Passing to a suitable version of $X$ if necessary to ensure
measurability of $d\Psi_x(X)$, we then have, in $([\mathbb{E}], e)$,
\[
\bigl[r_n \bigl( \Psi(Y_n) - \Psi(x) \bigr)\bigr]
\rightsquigarrow\bigl[d\Psi_x(X)\bigr], \qquad n \to\infty.
\]
\end{theorem}

%ad8.8 #&#
%
\begin{addendum}
\label{adddelta}
There exists a subset $\mathbb{W}_\infty$ of $\mathbb{W}$ and a
version $X'$ of $X$ such that the restriction of $d\Psi_x$ to $\mathbb
{W}_\infty$ is continuous, $X'$ takes values in $\mathbb{W}_\infty$
only, and $d\Psi_x(X')$ is Borel measurable.\vadjust{\goodbreak}
\end{addendum}

%%%%%%%%%%%%%%%%%%%%%%%%%%%%%%%%%%%%%%%%
%%% APPENDIX B
%%%%%%%%%%%%%%%%%%%%%%%%%%%%%%%%%%%%%%%%

%s9 #&#
\section{Proofs}
This appendix contains the most important proofs for the main part of
the paper, namely those for Sections~\ref{seccop}~and~\ref{secreg}. The remaining proofs are collected in Appendix~F
of the supplement [\citet{BucSegVol14supp}].

%%%%%%%%%%%%%%%%%%%%%%%%%%%%%%%%%%%%%%%%
%%% COPULAS PROOFS
%%%%%%%%%%%%%%%%%%%%%%%%%%%%%%%%%%%%%%%%
%s9.3 #&#
\subsection{Proofs for Section~\texorpdfstring{\protect\ref{seccop}}{4}}\label{appcop}
\mbox{}

\begin{pf*}{Proof of Theorem \ref{theocopweak}}
For the sake of a clear exposition, we split the proof into two
propositions. First, Proposition~\ref{prophypi} shows that
Condition~\ref{condD} implies a certain abstract hypi-differentia\-bility property stated in Condition~\ref{condhypidiff}. Then,
Proposition~\ref{propcopweakabs} shows that the latter condition
suffices to obtain the completion of Theorem~\ref{theocopweak}.
\end{pf*}

%co9.2 #&#
%
\begin{condition}[(Hypi-differentiability of $C$)] \label{condhypidiff}
Define the set
\[
\mathcal{W}(t):= \bigl\{a \in\bigl\{ \ell^\infty\bigl([0,1]\bigr) \bigr
\}^d\dvtx  \mathbf u + t a(\mathbf u) \in[0,1]^d\ \forall\mathbf u
\in[0,1]^d \bigr\},
\]
where $a(\mathbf u) = (a_1(u_1),\ldots, a_d(u_d))$.
Whenever $t_n \searrow0$, $t_n\ne0$, and $a_n=\break (a_{n1}, \ldots,
a_{nd}) \in\{\ell^\infty([0,1])\}^d$ converges uniformly to $a\in
\mathcal{W}:=\{ \mathcal{C}([0,1]) \}^d$ (i.e., $\| a_{nj}-a_j\|
_\infty\to0$ for all $j=1,\ldots, d$) such that $a_n \in\mathcal
{W}(t_n)$ for all $n\in\mathbb{N}$, the functions
\[
[0,1]^d \to\mathbb{R}\dvtx  u \mapsto% \Delta_{t_n}C(\vect u)(a_n(\vect u))=
t_n^{-1}
\bigl\{ C\bigl( \mathbf u + t_n a_n(\mathbf u) \bigr) - C(\mathbf u)
\bigr\}
\]
converge in $(\ell^\infty([0,1]^d), {{d_{\mathrm{hypi}}} })$ to
some limit ${d}C_a$.
% [0,1]^d \times\Wc\to\RR: (u,w) \mapsto\diff C_w(u).
\end{condition}

%pr9.3 #&#
%
\begin{proposition} \label{prophypi}
A copula $C$ satisfying Condition \ref{condD} also satisfies the
hypi-differentiability Condition~\ref{condhypidiff} with derivative
\[
{d}C_a(\mathbf u ) =\sup_{\varepsilon>0} \inf\Biggl\{
\sum_{j=1}^d \dot{C}_j(\mathbf v)
a_j(v_j)\dvtx  \mathbf v \in\mathcal{S}, | \mathbf{v} - \mathbf{u} | <
\varepsilon\Biggr\}, \qquad\mathbf{u} \in[0, 1]^d.
\]
\end{proposition}

Conversely, it is an open problem whether there exists a copula that
satisfies Condition~\ref{condhypidiff} but violates Condition~\ref
{condD}. According to the next proposition, Condition~\ref
{condhypidiff} can replace Condition~\ref{condD} in Theorem~\ref
{theocopweak}.

%pr9.4 #&#
%
\begin{proposition} \label{propcopweakabs}
Suppose that Condition~\ref{condweak} holds and that $C$ satisfies
Condition~\ref{condhypidiff}. Then, in $(L^\infty([0,1]^d),
{{d_{\mathrm{hypi}}} })$,
\[
\mathbb{C}_n \rightsquigarrow\mathbb{C}= \alpha+
{d}C_{(-\alpha_{1},\ldots,-\alpha_{d})}. %\quad\mbox{ in }~ (L^
\]
\end{proposition}

\begin{pf*}{Proof of Proposition \ref{prophypi}}
Let $t_n \searrow0$ and let $a_n \in\mathcal{W}(t_n)$ converge
uniformly to $a \in\mathcal{W}$. As in Condition~\ref
{condhypidiff}, we use the notation $a_n(\mathbf u) = (a_{n1}(u_1),\ldots, a_{nd}(u_d))$. We have to prove epi- and hypo-convergence of
\[
\mathbf u \mapsto F_n(\mathbf u):= t_n^{-1} \bigl\{ C
\bigl(\mathbf u+t_n a_n(\mathbf u)\bigr) - C(\mathbf u) \bigr\}
\]
to $F_\wedge$ and $F_\vee$, respectively, where $F = {d}C_a$.
Note that, in the notation of Appendix~\ref{appext}, we have
$F=G_\wedge^{\mathcal{S}\dvtx [0,1]^d}$, where
$G\dvtx \mathcal{S}\to\mathbb{R}$ is defined through
\[
G(\mathbf u) = \sum_{j=1}^d \dot
C_j(\mathbf u) a_j(u_j).
\]
By an application of Corollary~\ref{corextension}, and since
$F_\wedge=G_\wedge^{\mathcal{S}\dvtx [0,1]^d}$ and
$F_\vee=G_\vee^{\mathcal{S}\dvtx [0,1]^d}$,
it suffices to show that:
\begin{longlist}[(ii)]
\item[(i)] $\forall\mathbf u \in[0,1]^d\dvtx  \forall\mathbf u_n \to\mathbf u\dvtx
\liminf_{n\to\infty} F_n(\mathbf u_n) \ge F_\wedge(\mathbf u)$,
\item[(ii)] $\forall\mathbf u \in[0,1]^d\dvtx  \forall\mathbf u_n \to\mathbf u\dvtx
\limsup_{n\to\infty} F_n(\mathbf u_n) \le F_\vee(\mathbf u)$.
\end{longlist}

We begin with the proof of (i) and fix a point $\mathbf u\in[0,1]^d$ and a
sequence $\mathbf u_n\to\mathbf u$. Choose $\varepsilon> 0$ and let $| \cdot
|_1$ denote the $L_1$-norm on $\mathbb{R}^d$.
Due to Lemma \ref{lemcav}, we may choose
\[
\mathbf u_n^\star\in\bigl\{ \mathbf v\in[0,1]^d\dvtx  | \mathbf
u_n - \mathbf v|_1 \le\varepsilon t_n/2 \bigr\}
\]
and
\[
\mathbf u_n^\circ\in\bigl\{ \mathbf v\in[0,1]^d\dvtx  \bigl| \mathbf
u_n + t_n a_n(\mathbf u_n) - \mathbf{v}
\bigr|_1 \le\varepsilon t_n/2\bigr\}
\]
such that, for the path
\[
\gamma_n(s)=(1-s) \mathbf u_n^\star+ s \mathbf
u_n^\circ, \qquad s\in[0,1],
\]
the set $\{s \in[0,1]\dvtx  \gamma_n(s) \notin\mathcal{S}\}$ has
Lebesgue-measure zero. Define
$
f_n(s) = t_n^{-1} C(\gamma_n(s))$, $s \in[0, 1]$,
and note that
\begin{eqnarray*}
\bigl\llvert\bigl\{f_n(1) - f_n(0) \bigr\} -
F_n(\mathbf u_n) \bigr\rrvert
&=& t_n^{-1} \bigl\llvert C\bigl(\mathbf u_n^\circ
\bigr) - C\bigl(\mathbf u_n+t_na_n(\mathbf
u_n)\bigr) - C\bigl(\mathbf u_n^\star\bigr) + C(\mathbf
u_n) \bigr\rrvert
\\
&\le& t_n^{-1} \bigl\{ \bigl|\mathbf u_n +
t_n a_n(\mathbf u_n) - \mathbf u_n^\circ\bigr|_1
+ \bigl| \mathbf u_n^\star- \mathbf u_n\bigr|_1 \bigr
\} \le\varepsilon
\end{eqnarray*}
by Lipschitz-continuity of $C$.
Lipschitz-continuity of $C$ also implies absolute continuity of $f_n$, which
allows us to choose $\mathbf v_n\in\gamma_n([0,1])\cap\mathcal{S}$ such that
% \sv{strictly formally, the derivative of the copula needs not be
%defined outside of $\mathcal{S}$. It is clear what we mean when we
%write the integral, but perhaps we should still comment on that?}
%
\begin{eqnarray*}
\varepsilon+ F_n(u_n)
&\ge& f_n(1)- f_n(0)
\\
& =& \int_{0}^{1}
f_n'(s) \,ds %\\
\\
&=& \sum
_{j=1}^d t_n^{-1}
\bigl(u_{nj}^\circ- u_{nj}^\star\bigr) \int
_0^1 \dot C_j\bigl(
\gamma_n(s)\bigr) \,ds
\\
&=& \sum_{j=1}^d \bigl[
a_{nj}(u_{nj}) + t_n^{-1} \bigl\{
u_{nj}^\circ- u_{nj}^\star- t_n
a_{nj}(u_{nj}) \bigr\} \bigr] \int_0^1
\dot C_j\bigl(\gamma_n(s)\bigr) \,ds
\\
& \ge&\inf_{s\dvtx  \gamma_n (s) \in\mathcal{S}} \sum_{j=1}^d
a_{nj}(u_{nj}) \dot C_j\bigl(
\gamma_n(s)\bigr) - \varepsilon
\\
& \ge&\sum_{j=1}^d a_{nj}(v_{nj})
\dot C_j(\mathbf v_n) + \sum_{j=1}^{d}
\bigl\{ a_{nj}(u_{nj}) - a_{nj}(v_{nj})
\bigr\} \dot C_j(\mathbf v_n) - 2\varepsilon
\\
& \ge&\sum_{j=1}^d a_{nj}(v_{nj})
\dot C_j(\mathbf v_n) - 3\varepsilon=F(\mathbf v_n)
- 3\varepsilon= F_\wedge(\mathbf v_n) - 3 \varepsilon
\end{eqnarray*}
for sufficiently large $n$, where we have used the bounds $0 \le\dot
C_j \le1$, uniform convergence of $a_{nj}$ to $a_j$, uniform
continuity of~$a_j$ and the fact that $F$ is continuous in $\mathbf v_n$.
Hence, by lower semicontinuity of $F_\wedge$,
\[
\liminf_{n\to\infty} F_n(\mathbf u_n) \ge
F_\wedge(\mathbf u) - 4\varepsilon.
\]
As $\varepsilon>0$ was arbitrary the assertion in (i) follows.

The proof of (ii) is analogous. In the main inequality chain, all signs
can be reversed if the infimum is replaced by a supremum and upon
noting that on $\mathcal{S}$, the functions $F$, $F_\wedge$ and
$F_\vee$ are equal.
%Lipschitz-continuity of $C$ also implies absolute continuity of the
%function $f_n$
%whence we can write
%%\sv{strictly formally, the derivative of the copula needs not be
%defined outside of $\mathcal{S}$. It is clear what we mean when we
%write the integral, but perhaps we should still comment on that?}
% f_n(1)- f_n(0) = \int_{0}^{1} f_n'(s) \,ds = \sum_{j=1}^d
%t_n^{-1}(u_{nj}^\circ- u_{nj}^\star) \int_0^1 \dot C_j(\gamma_n(s))
% = \sum_{j=1}^d \left[ a_{nj}(u_{nj}) + t_n^{-1} \left\{u_{nj}^\circ-
%u_{nj}^\star- t_n a_{nj}(u_{nj}) \right\}\right] \int_0^1 \dot C_j(
%Thus, we can choose some $\vect v_n\in\gamma_n([0,1])\cap\Sc$ such
%that the right-hand side of the previous expression is bounded from
%below by
% f_n(1)- f_n(0)
% &\ge\inf_{s\dvtx  \gamma_n (s) \in\Sc} \sum_{j=1}^d a_{nj}(u_{nj}) \dot
%C_j(\gamma_n(s)) - \eps\\
% & \ge&\sum_{j=1}^d a_{nj}(v_{nj}) \dot C_j(\vect v_n) +
%v_n) - 2\eps\\
% & \ge&\sum_{j=1}^d a_{nj}(v_{nj}) \dot C_j(\vect v_n) - 3\eps=F(
%for sufficiently large $n$, by boundedness of $\dot C_j$ as a
%consequence of the Lipschitz-continuity of $C$, uniform convergence of
%$a_{nj}$ and uniform continuity of $a_j$.
%Hence, by lower semi-continuity of $F_\wedge$,
% \liminf_{n\to\infty} F_n(\vect u_n) \ge F_\wedge(\vect u) - 4\eps
%and as $\eps>0$ was arbitrary the assertion in (i) follows.
%The proof of (ii) is analogous. \ab{Last step in inequality chain.}
\end{pf*}

%For $0<\eps<1/2$ let $\vect u, \vect v\in[\eps,1-\eps]^d$ be two
%distinct points and denote by $H_{\vect u}$ and $H_{\vect v}$ the
%hyperplanes being orthogonal to $\vect u- \vect v$ and passing through
%$\vect u$ and $\vect v$, respectively. For $0<\delta<\eps$ define $H_{
%v}^\delta= H_{\vect v}\cap B_1(\vect v,\delta)$, where $B_1(\vect u,
%u$ with respect to the $\| \cdot\|_1$-norm. Finally, let $\Zc$ denote
%the cylinder of height $h = \| \vect u- \vect v\|$ and radius $\delta$
%with area of the top equal to $H_{\vect u}^\delta$ and area of the
%bottom equal to $H_{\vect v}^\delta$, i.e.,
% \Zc= \{ \vect x \in\RR^d\dvtx  \vect x = \vect y + s(\vect v-\vect u),
%Let $\Dc$ be a Lebesgue-null set in $[0,1]^d$ and define, for any $
% \Zc_{\vect y}^\Dc=\left\{ s \in\RR\dvtx  \vect y + s(\vect v- \vect u)
%Then $\Zc_{\vect y}^\Dc$ is a one-dimensional Lebesgue-null set for
%almost all $\vect y\in H_{\vect u}^\delta$.
%
%By translation, it suffices to prove the result for the standard
%cylinder, i.e., for $\vect u=(0,\ldots, 0)$ and $\vect v=(0,
% \Zc= \{ \vect x\in\RR^{d-1}\dvtx  \|\vect x\|_1 \le r\} \times[0,h].
%By Cavalieri's principle,
% \lambda_d(\Zc\cap\Dc) = \int_{\RR^{d-1}} \lambda_1(\{ s \in\RR\dvtx  (
% = \int_{\RR^{d-1}} \lambda_1(\Zc_{(\vect x,0)}^\Dc) \,d\vect x
%and since the expression on the left is equal to $0$ the assertion
%follows.

\begin{pf*}{Proof of Proposition~\ref{propcopweakabs}}
Recall %the definition of
$\beta_n=(\beta_{n1},\ldots, \beta_{nj})$, with $\beta_{nj} =
\sqrt n ( G_{nj}^- - \mathrm{id}_{[0,1]})$. It follows from
Condition~\ref{condweak} and the functional delta method for the
inverse mapping, also known as Vervaat's lemma, that
\[
(\alpha_n, \beta_n) = (\alpha_n,
\beta_{n1},\ldots, \beta_{nd}) \rightsquigarrow(\alpha, -
\alpha_1,\ldots, -\alpha_d)
\]
in $\ell^\infty([0,1]^d) \times\{ \ell^\infty([0,1]) \}^d$, with
respect to the supremum distance in each coordinate. Note that we can write
$\mathbb{C}_n=g_n(\alpha_n, \beta_n)$, where $g_n\dvtx  \ell^\infty
([0,1]^d) \times\mathcal{W}(1/\sqrt n) \to(L^\infty([0,1]^d),
{{d_{\mathrm{hypi}}} })$ is defined as
%e9.2 #&#
%
\begin{equation}
\label{eqgn} g_n(a,b) = a(\mathrm{id}_{[0,1]^d} + b/\sqrt n) +
\sqrt n \bigl\{ C(\mathrm{id}_{[0,1]^d} +b/\sqrt n) - C\bigr\}.
\end{equation}
Exploiting Condition \ref{condhypidiff} and Lemma~\ref{lemsum}
(recall that $\alpha$ is continuous almost surely), the assertion
follows from the extended continuous mapping theorem, see
Theorem~1.11.1 in \citet{vandwell1996}.
\end{pf*}

%le9.5 #&#
%
\begin{lemma} \label{lemcav}
Let $\mathbf u, \mathbf v \in\mathbb{R}^d$ be two distinct points and denote
by $H_{\mathbf u}$ and $H_{\mathbf v}$ the hyperplanes being orthogonal to $\mathbf
u- \mathbf v$ and passing through $\mathbf u$ and $\mathbf v$, respectively. For
$\delta> 0$, set $H_{\mathbf u}^\delta= H_{\mathbf u} \cap B_1(\mathbf u, \delta
)$ and $H_{\mathbf v}^\delta= H_{\mathbf v}\cap B_1(\mathbf v,\delta)$, where
$B_1(\mathbf u,\delta)$ denotes the unit ball of radius $\delta$ centered
at $\mathbf u$ with respect to the $\| \cdot\|_1$-norm. Finally, let
$\mathcal{Z}$ denote the cylinder
% of height $h = \| \vect u- \vect v\|$ and radius $\delta$
with top area equal to $H_{\mathbf u}^\delta$ and bottom area equal to
$H_{\mathbf v}^\delta$, that is,
\[
\mathcal{Z}= \bigl\{ \mathbf y + s(\mathbf v-\mathbf u)\dvtx  \mathbf y \in H_{\mathbf u}^\delta,
s \in[0,1] \bigr\}.
\]
Let $\mathcal{D}$ be a Lebesgue-null set in $\mathbb{R}^d$ and
define, for any $\mathbf y \in H_{\mathbf u}^\delta$,
\[
\mathcal{Z}_{\mathbf y}^\mathcal{D}= \bigl\{ s \in\mathbb{R}\dvtx  \mathbf y +
s(\mathbf v- \mathbf u) \in\mathcal{Z}\cap\mathcal{D} \bigr\}.
\]
Then $\mathcal{Z}_{\mathbf y}^\mathcal{D}$ is a one-dimensional
Lebesgue-null set for almost all $\mathbf y\in H_{\mathbf u}^\delta$.
\end{lemma}

The proofs of Lemma~\ref{lemcav}, Propositions~\ref{propbootweak}~and~\ref{proploc} are given in Appendix~F.3
in the supplement [\citet{BucSegVol14supp}].

\subsection{Proofs for Section~\texorpdfstring{\protect\ref{secreg}}{6}}\label{appreg}
\mbox{}

\begin{pf*}{Proof of Theorem \ref{theoresweak}}
The proof consists of two main steps. In the first step, consider $\ell
^\infty( \bar\mathbb{R}) \times\mathbb{R}^p$ equipped with the metric
\[
\rho\bigl( (h_1, \mathbf{y}_1), (h_2,
\mathbf{y}_2) \bigr) = \| h_1 - h_2
\|_\infty+ | \mathbf{y}_1 - \mathbf{y}_2 |.
\]
As shown in Appendix~F.5 in the supplement [\citet
{BucSegVol14supp}], we have, in $(\ell^\infty( \bar\mathbb{R})
\times\mathbb{R}^p, \rho)$, as $n\to\infty$,
%e9.4 #&#
%
\begin{eqnarray}\label{eqdweak}
\qquad \bigl( \mathbb{G}_n f_{\cdot, \hat{\bolds{\beta}}_n - \bolds{\beta}}, \sqrt
{n} ( \hat{\bolds{
\beta}}_n - \bolds{\beta} ) \bigr) = ( \mathbb{G}_n
f_{\cdot, \mathbf{0}}, \mathbb{G}_n \bolds{\psi} ) + o_p(1)
&\rightsquigarrow &( \mathbb{G}f_{\cdot, \mathbf{0}},
\mathbb{G}\bolds{
\psi} )
\end{eqnarray}
%
% \bigl(
% \GG_n f_{\cdot, \hat{\vect{\beta}}_n - \vect{\beta}},
% \sqrt{n} ( \hat{\vect{\beta}}_n - \vect{\beta} )
% \bigr)
% &=&
% \bigl(
% \GG_n f_{\cdot, \vect{0}},
% \GG_n \vect{\psi}
% \bigr)
% +
% o_p(1) \\
% & \weak
% \bigl(
% \GG f_{\cdot, \vect{0}},
% \GG\vect{\psi}
% \bigr),
% \qquad n \to\infty
where $\mathbb{G}$ denotes a zero-mean Gaussian process on $\mathcal
{G}= \mathcal{F}\cup\{ \psi_1, \ldots, \psi_p \} $ with covariance
given in (\ref{eqcovG}).
%The first coordinate of the limit, $\GG f_{\cdot, \vect{0}}$, is
%continuous in $z$ almost surely by Lemma~\ref{lemcont} below.
Define $T_n\dvtx  \ell^\infty( \bar\mathbb{R}) \times\mathbb{R}^p \to
\ell^\infty( \bar\mathbb{R})$ by
\[
T_n( G, \bolds{\gamma} ) = G + g_n(\bolds{\gamma}),
\]
where the map $g_n(\gamma)\in\ell^\infty(\bar\mathbb{R})$ is
defined by $(g_n(\gamma))(\pm\infty) = 0$ and
\begin{eqnarray*}
\bigl(g_n(\bolds{\gamma})\bigr) (z) &=& t_n^{-1}
\int_{\mathbb{R}^p} \bigl\{ F\bigl(z + t_n
\mathbf{x}' \bolds{\gamma}\bigr) - F(z) \bigr\} P^X( {d}
\mathbf{x} ), \qquad z \in\mathbb{R}.
\end{eqnarray*}
Note that we can write the second term in (\ref{eqFndecomp}) as
%e9.5 #&#
%
\begin{equation}
\label{eqgnWhy} \sqrt{n} \{ P f_{z, \hat{\bolds{\beta}}_n-\bolds\beta} - P
f_{z, \mathbf
{0}} \} =
\bigl(g_n \bigl( \sqrt{n} (\hat{\bolds{\beta}}_n - \bolds{
\beta}) \bigr) \bigr) (z)
\end{equation}
with $t_n = 1/\sqrt{n}$. This also allows to write
\[
\mathbb{F}_n = T_n \bigl( \mathbb{G}_n
f_{\cdot, \hat{\bolds{\beta}}_n - \bolds{\beta}}, \sqrt{n} ( \hat{\bolds{\beta
}}_n - \bolds{\beta} )
\bigr).
\]
The assertion of Theorem~\ref{theoresweak} will then follow by an
application of the extended continuous mapping theorem.
More precisely, if $G_n, G \in\ell^\infty( \bar\mathbb{R})$ are
such that $G$ is continuous and $\| G_n - G \|_\infty\to0$, and if
moreover $\bolds{\gamma}_n \to\bolds{\gamma}$ in $\mathbb{R}^p$, then,
in $(\ell^\infty( \bar\mathbb{R}), {d_{\mathrm{hypi}}} )$,
%e9.6 #&#
%
\begin{equation}
\label{eqcontco} T_n( G_n, \bolds{\gamma}_n ) \to
T( G, \bolds{\gamma} ):= G + g(\bolds{\gamma}),
\end{equation}
by Lemma~\ref{lemgng} below and Lemma~\ref{lemsum} on weak
hypi-convergence of sums.
Here, the map $g(\bolds\gamma) \in\ell^\infty(\bar\mathbb{R})$ is
defined by $g(\bolds\gamma)(\pm\infty)=0$ and, for $z\in\mathbb{R}$,
%e9.7 #&#
%
\begin{eqnarray}\label{eqg}
\bigl(g(\bolds{\gamma})\bigr) (z) &=& - f(z-) \int
_{-\infty}^0 \mathbb{P}\bigl( \mathbf{X}'
\bolds{\gamma} < y \bigr) \,{d}y
\nonumber\\[-8pt]\\[-8pt]
&&{} + f(z+) \int_0^{+\infty}
\mathbb{P}\bigl( \mathbf{X}' \bolds{\gamma} > y \bigr)\, {d}y.\nonumber
%, \qquad z \in\reals.
\end{eqnarray}
%
% (g(\vect{\gamma}))(z)
% =
% -
% f(z-) \int_{-\infty}^0 \Pr( \vect{X}' \vect{\gamma} < y ) \diff y % %
%+
% f(z+) \int_0^{+\infty} \Pr( \vect{X}' \vect{\gamma} > y ) \diff y
%, \qquad z \in\reals.
Note that the integrals on the right-hand side of the last display
exist as a consequence of condition~(R3) and Fubini's theorem, which
implies that
%e9.8 #&#
%e9.9 #&#
%
\begin{eqnarray}
\label{eqintegralplus} \int_0^{+\infty} \mathbb{P}\bigl(
\mathbf{X}' \bolds{\gamma} > y \bigr) \,{d}y &=& \mathbb{E}\bigl[ \max
\bigl( \mathbf{X}' \bolds{\gamma}, 0 \bigr) \bigr]<\infty,
\\
\label{eqintegralminus} \int_{-\infty}^0 \mathbb{P}\bigl(
\mathbf{X}' \bolds{\gamma} < y \bigr) \,{d}y &=& \mathbb{E}\bigl[ \max
\bigl( - \mathbf{X}' \bolds{\gamma}, 0 \bigr) \bigr]<\infty.
\end{eqnarray}

Finally, as a consequence of (\ref{eqcontco}) and since $\mathbb
{G}f_{\cdot, \mathbf{0}}$ is continuous almost surely by Lemma~F.5 in the supplementary material, the assertion follows from
(\ref{eqdweak}) and an application of the extended continuous mapping
theorem [\citet{vandwell1996}, Theorem~1.11.1].
\end{pf*}

%.; its proof is given in Appendix~\ref{appreg2} in the supplement [

The preceding proof made use of Lemma~\ref{lemgng} below. For its
formulation, we need two additional lemmas. The proof of the first one
is trivial and, therefore, omitted.

%First we consider the term $\GG_n f_{z, \hat{\vect{\beta}}_n - \vect{
% \Fclass
% = \{ f_{z, \vect{\delta}}\dvtx  z \in\reals, \vect{\delta} \in
%Clearly, $\Fclass\subset L^2(P)$.

%le9.7 #&#
%
\begin{lemma}
\label{lemladlag}
If $f$ is l\`adl\`ag, then both functions $z \mapsto f(z+)$ and $z
\mapsto f(z-)$ defined in \textup{(R2)} of Theorem~\ref{theoresweak} are l\`
adl\`ag, too. Their right-hand limits at $z$ are both equal to $f(z+)$
and their left-hand limits at $z$ are both equal to $f(z-)$.
% \lim_{0 < s \to0} f((z + s)\pm)
%% &=& \lim_{0 < s \to0} f((z + s)-)
% = f(z+),
% \lim_{0 < s \to0} f((z - s)\pm)
%% &=& \lim_{0 < s \to0} f((z - s)-)
% = f(z-).
% \lim_{0 < s \to0} f((z + s)+)
% &=& \lim_{0 < s \to0} f((z + s)-)
% = f(z+), \\
% \lim_{0 < s \to0} f((z - s)+)
% &=& \lim_{0 < s \to0} f((z - s)-)
% = f(z-).
\end{lemma}

%le9.8 #&#
%
\begin{lemma}
\label{eqgupperlower}
%If Conditions~\ref{condladlag} and~\ref{condXL1} hold,
If conditions~\textup{(R2)} and \textup{(R3)} hold,
then for every $\gamma\in\mathbb{R}^p$, the function $g(\bolds{\gamma
})$ in~(\ref{eqg}) is uniformly bounded and l\`adl\`ag, with right-
and left-hand limits at $z \in\mathbb{R}$ given by
%e9.10 #&#
%
\begin{eqnarray}
\label{eqggammaplusminus} \bigl(g(\bolds{\gamma})\bigr) (z\pm) &=& f(z\pm
) \mathbb{E}\bigl[
\mathbf{X}' \bolds{\gamma} \bigr].
\end{eqnarray}
%
% (g(\vect{\gamma}))(z+)
% &=& f(z+) \expec[ \vect{X}' \vect{\gamma} ], \\
% (g(\vect{\gamma}))(z-)
% &=& f(z-) \expec[ \vect{X}' \vect{\gamma} ].
The upper and lower semicontinuous hulls of $g(\bolds{\gamma})$ at $z
\in\mathbb{R}$ are
%e9.11 #&#
%e9.12 #&#
%
\begin{eqnarray}
\label{eqggammaupper} \bigl(g(\bolds{\gamma})\bigr)_{\vee}(z) &=& \max
\bigl\{
\bigl(g(\bolds{\gamma})\bigr) (z-), \bigl(g(\bolds{\gamma})\bigr) (z), \bigl
(g(\bolds{
\gamma})\bigr) (z+) \bigr\},
\\
\label{eqggammalower} \bigl(g(\bolds{\gamma})\bigr)_{\wedge}(z) &=& \min
\bigl\{
\bigl(g(\bolds{\gamma})\bigr) (z-), \bigl(g(\bolds{\gamma})\bigr) (z), \bigl
(g(\bolds{
\gamma})\bigr) (z+) \bigr\}.
\end{eqnarray}
Moreover, $(g(\bolds{\gamma}))_{\wedge}(\pm\infty) = (g(\bolds{\gamma
}))_{\vee}(\pm\infty) = 0$.
\end{lemma}

\begin{pf} %{Proof of Lemma~\ref{eqgupperlower}}
The existence and the expressions of the right-hand and left-hand
limits of $g(\bolds{\gamma})$ at $z \in\mathbb{R}$ are a consequence
of Lemma~\ref{lemladlag} and the fact that
\[
- \int_{-\infty}^0 \mathbb{P}\bigl(
\mathbf{X}' \bolds{\gamma} < y \bigr) \,{d}y + \int
_0^{+\infty} \mathbb{P}\bigl( \mathbf{X}'
\bolds{\gamma} > y \bigr) \,{d}y = \mathbb{E}\bigl[ \mathbf{X}' \bolds{
\gamma} \bigr],
\]
which follows in turn from (\ref{eqintegralplus}) and (\ref
{eqintegralminus}). The statement about the upper (lower)
semicontinuous hull follows from the fact that for a l\`adl\`ag
function, the supremum (infimum) over a shrinking neighborhood around a
point converges to the maximum (minimum) of the function value at the
point itself and the right-hand and left-hand limits at that point.
\end{pf}
%
%}

%le9.9 #&#
%
\begin{lemma}
\label{lemgng}
%Assume Conditions~\ref{condladlag} and~\ref{condXL1}.
Assume conditions~\textup{(R2)} and~\textup{(R3)} in Theorem~\ref{theoresweak}.
If $\bolds{\gamma}_n \to\bolds{\gamma}$ in~$\mathbb{R}^p$, then
\[
{d_{\mathrm{hypi}}} \bigl( g_n(\bolds{\gamma}_n), g(\bolds{
\gamma}) \bigr) \to0, \qquad n \to\infty.
\]
%
%where $\dhypi$ denotes any semimetric inducing hypi-convergence on $
\end{lemma}

\begin{pf}
First of all, for $z \in\mathbb{R}$, we can write $(g_n(\bolds{\gamma
}))(z)$ as
\[
\int_{\mathbb{R}^p} t_n^{-1}
\int_0^{t_n \mathbf{x}' \bolds{\gamma}} f(z + y) \,{d}y P^X(
{d}\mathbf{x} ) = \int_{\mathbb{R}^p} \int_0^{\mathbf{x}' \bolds{\gamma}}
f(z + t_n y) \,{d}y P^X({d}\mathbf{x}).
\]
It follows that $g_n(\bolds{\gamma}_n)$ is uniformly close to $g_n(\bolds
{\gamma})$:
%letting $\| \cdot\|_\infty$ denote the supremum norm for
%real-valued functions on $\reals$,
we have
\[
\bigl| \bigl(g_n(\bolds{\gamma}_n)\bigr) (z) -
\bigl(g_n(\bolds{\gamma})\bigr) (z) \bigr| \le\| f \|_\infty\int
_{\mathbb{R}^p} | \mathbf{x} | P^X( {d}\mathbf{x} ) | \bolds{
\gamma}_n - \bolds{\gamma} |, \qquad z \in\mathbb{R},
\]
and thus, noting that $(g_n(\bolds{\gamma}_n))(\pm\infty) = 0 =
(g_n(\bolds{\gamma}))(\pm\infty)$,
%e9.13 #&#
%
\begin{equation}
\label{eqgngammangngamma} \bigl\| g_n(\bolds{\gamma}_n) - g_n(
\bolds{\gamma}) \bigr\|_\infty\to0, \qquad n \to\infty.
\end{equation}
Hence, without loss of generality, we can assume that $\bolds{\gamma}_n
= \bolds{\gamma}$ for all $n$.

Fix $z \in\bar\mathbb{R}$. We will prove hypi-convergence of
$g_n(\bolds{\gamma})$ to $g(\bolds{\gamma})$
%by proving hypo-convergence of $g_n(\vect{\gamma})$
%to $(g(\vect{\gamma}))_{\vee}$ at $z$ and epi-convergence of $g_n(
%of~\eqref{eqdefepi} and \eqref{eqdefhypo}.
by using the pointwise criteria in~(\ref{eqdefepi}) and (\ref{eqdefhypo}).
First, consider the case $z_n \to z = +\infty$. %By Lemma~
Observe that for any fixed $\mathbf{x} \in\mathbb{R}^p$
\[
t_n^{-1}\bigl|F\bigl(z_n+t_n
\mathbf{x}'\bolds{\gamma}\bigr) - F(z_n)\bigr| \leq\bigl|
\mathbf{x}'\bolds{\gamma}\bigr| \sup_{y\geq z_n-t_n|\mathbf{x}'\bolds{\gamma}|} f(y) \to0
\]
%
%since the existence of the limit $\lim_{z\to+\infty}f(z)$ implies that
%this limit must be zero [recall that $f$ is a density] and thus $
since $\lim_{z\to+\infty}f(z)=0$ by assumption.
Hence, by dominated convergence, $(g_n(\bolds{\gamma}))(z_n) \to0$,
which is equal to $(g(\bolds{\gamma}))_{\wedge}(+\infty) = (g(\bolds
{\gamma}))_{\vee}(+\infty)$ in view of Lemma~\ref{eqgupperlower}.
The limit $z_n \to-\infty$ can be handled similarly.

It thus remains to consider $z_n \to z \in\mathbb{R}$. By Fubini's
theorem, we have
\begin{eqnarray*}
\bigl(g_n(\bolds{\gamma})\bigr) (z_n)
&=& - \int\!\!\int_{\mathbf{x}' \bolds{\gamma} < y < 0} f(z_n + t_n
y) \,{d}y P^X( {d}\mathbf{x} )
\\
&&{} + \int\!\!\int
_{\mathbf{x}' \bolds{\gamma} > y > 0} f(z_n + t_n y) \,{d}y
P^X( {d}\mathbf{x} )
\\
&=& - \int_{-\infty}^0 f(z_n +
t_n y) \mathbb{P}\bigl( \mathbf{X}' \bolds{\gamma} < y \bigr)\,
{d}y
\\
&&{} + \int_0^{+\infty} f(z_n +
t_n y) \mathbb{P}\bigl( \mathbf{X}' \bolds{\gamma} > y \bigr)
\,{d}y.
\end{eqnarray*}
The idea is now to replace $f(z_n + t_n y)$ by $f(z-)$ or $f(z+)$
according to whether $z_n + t_n y$ is smaller or larger than $z$. To
this end, define the auxiliary functions
\begin{eqnarray*}
w(y) &=& \cases{ \displaystyle- \mathbb{P}\bigl( \mathbf{X}' \bolds{
\gamma} < y \bigr), &\quad if $y < 0$,
\vspace*{5pt}\cr
\displaystyle\mathbb{P}\bigl(
\mathbf{X}' \bolds{\gamma} > y\bigr), &\quad if $y > 0$;}
\\
\eta(a) &=& f(z-) \int_{-\infty}^a w(y) \,{d}y
+ f(z+) \int_a^{+\infty} w(y) \,{d}y, \qquad a
\in\mathbb{R}.
\end{eqnarray*}
Further, put $a_n= (z-z_n)/t_n$
% a_n = \frac{z - z_n}{t_n}
and observe that $z_n + t_n y < z$ if $y < a_n$ while $z_n + t_n y > z$
if $y > a_n$. We have
\[
\bigl(g_n(\bolds{\gamma})\bigr) (z_n) = \int
_{-\infty}^\infty f(z_n + t_n y)
w(y) \,{d}y.
\]
By the dominated convergence theorem, as $\int_\mathbb{R}|w(y)|
\,{d}y = \mathbb{E}[ | \mathbf{X}' \bolds{\gamma} | ] < \infty$,
%e9.14 #&#
%
\begin{eqnarray}\label{eqgnEtaToZero}
&& \bigl(g_n(\bolds{\gamma})\bigr) (z_n) -
\eta(a_n)\nonumber
\\
&&\qquad = \int_{-\infty}^{+\infty} \bigl\{
f(z_n + t_n y) - f(z-) \bigr\} \mathbh{1}( y <
a_n ) w(y) \,{d}y
\\
&&\quad\qquad{}  + \int_{-\infty}^{+\infty}
\bigl\{ f(z_n + t_n y) - f(z+) \bigr\} \mathbh{1}( y >
a_n ) w(y) \,{d}y =o(1).\nonumber
\end{eqnarray}

Consider the extrema of the function $\eta$. The function $\eta$ can
be written as
\[
\eta(a) = f(z+) \int_{-\infty}^{+\infty} w(y) \,{d}y +
\bigl\{ f(z-) - f(z+) \bigr\} \int_{-\infty}^a w(y)
\,{d}y.
\]
It follows that $\eta$ is absolutely continuous with Radon--Nikodym
derivative % $\dot{\eta}$ given by
$
\dot{\eta}(a) = \{ f(z-) - f(z+) \} w(a)$.
Since $w(y) \le0$ for $y < 0$ and $w(y) \ge0$ for $y > 0$, we find
that $\eta$ is monotone on $(-\infty, 0)$ and on $(0, \infty)$.
Hence, $\eta$ attains its extrema at either $a \to-\infty$, $a = 0$,
or $a \to+\infty$. But for $a \to\pm\infty$, we find from (\ref
{eqggammaplusminus}) that
\[
\eta(\mp\infty) = f(z\pm) \int_{-\infty}^{+\infty} w(y)
\,{d}y = f(z\pm) \mathbb{E}\bigl[ \mathbf{X}' \bolds{\gamma} \bigr] =
\bigl(g(\bolds{\gamma})\bigr) (z\pm),
\]
while for $a = 0$, we find
\begin{eqnarray*}
\eta(0) &=& f(z-) \int_{-\infty}^0 w(y) \,{d}y +
f(z+) \int_0^{+\infty} w(y) \,{d}y = \bigl(g(\bolds{
\gamma})\bigr) (z).
\end{eqnarray*}
As a consequence, using (\ref{eqggammaupper}),
\begin{eqnarray*}
\sup_{a \in\mathbb{R}} \eta(a) &=& \max\bigl\{ \eta(-\infty), \eta(0),
\eta(+\infty) \bigr\}
\\
&=& \max\bigl\{ \bigl(g(\bolds{\gamma})\bigr) (z-), \bigl(g(\bolds{\gamma
})\bigr) (z),
\bigl(g(\bolds{\gamma})\bigr) (z+) \bigr\} = \bigl(g(\bolds{\gamma})
\bigr)_{\vee}(z),
\end{eqnarray*}
and similarly, by (\ref{eqggammalower}),
$
\inf_{a \in\mathbb{R}} \eta(a)
= (g(\bolds{\gamma}))_{\wedge}(z)$.
In combination with (\ref{eqgnEtaToZero}), we obtain that
\begin{eqnarray*}
\bigl(g(\bolds{\gamma})\bigr)_{\wedge}(z) &=& \inf_{a \in\mathbb{R}}
\eta(a)
\\
& \stackrel{(I1)} {\le}& \liminf_{n \to\infty}
\bigl(g_n(\bolds{\gamma})\bigr) (z_n) \le\limsup
_{n \to\infty} \bigl(g_n(\bolds{\gamma})\bigr)
(z_n)
\\
&\stackrel{(I2)} {\le}& \sup_{a \in\mathbb{R}} \eta(a) = \bigl(g(\bolds{
\gamma})\bigr)_{\vee}(z).
\end{eqnarray*}
Moreover, the inequalities $(I1)$ and $(I2)$ in the above display
become equalities if we choose $z_n$ in such a way that $a_n = (z -
z_n) / t_n$ converges to $-\infty$, $0$, or $\infty$, according to
where the infinimum and supremum of $\eta$ are attained.

The above paragraph shows that $g_n(\bolds{\gamma})$ epi-converges to
$(g(\bolds{\gamma}))_{\wedge}$ and hypo-converges to $(g(\bolds{\gamma
}))_{\vee}$ pointwise at every $z \in\mathbb{R}$. As a consequence,
$g_n(\bolds{\gamma})$ hypi-converges to $g(\bolds{\gamma})$. This
completes the proof.
%Since moreover the supremum distance between $g_n(\vect{\gamma}_n)$
%and $g_n(\vect{\gamma})$ converges to zero by
%hypi-converges to $g(\vect{\gamma})$. The proof of the Lemma is
%complete.
\end{pf}
\end{appendix}

% zodis "Acknowledgments" paliekamas pagal autoriu
\section*{Acknowledgments} The authors would like to thank the
Associate Editor and two referees for their constructive comments on an
earlier version of this manuscript.

Parts of this paper were written when A. B\"ucher was a post-doctoral
researcher at Universit\'e catholique de Louvain, Belgium.

\begin{supplement}[id=suppA]
% \sname{Proof of auxiliary results}
\stitle{Supplement to: ``When uniform weak convergence fails: Empirical
processes for dependence functions and residuals~via epi- and hypographs''}
\slink[doi]{10.1214/14-AOS1237SUPP} %[doi,text={...}] - jei reikia
%suskaldyti doi
\sdatatype{.pdf}
\sfilename{aos1237\_supp.pdf}
\sdescription{In the supplement, missing proofs for the results in
this paper are given.}
\end{supplement}

% imsref loaded by linak, 2014-06-20 15:06:08
%
% imsref loaded by linak, 2014-06-26 14:35:35

\printaddresses
\end{document}